\documentclass{amsart}
\usepackage{amsmath, amstext, amsxtra, amsthm, amscd, amsgen,amsbsy, amsopn, amsfonts, latexsym, amssymb, amscd}
\usepackage{marvosym}
\usepackage[all]{xy}

\usepackage{color}
\usepackage{multicol}
\usepackage{colortbl}

\begin{document}

%%%%%%%%%%%%%MACROS%%%%%%%%%%%%%%%%%%%%%%%%%%%%%%%%%%%%%%%%%%%%%

\def\wrt{with respect to}

\def\pproof{\noindent{\bf Proof:} }

\def\proofend{\hbox to 1em{\hss}\hfill $\blacksquare $\bigskip }

\newtheorem{theorem}{Theorem}[section]
\newtheorem{proposition}[theorem]{Proposition}
\newtheorem{lemma}[theorem]{Lemma}
\newtheorem{remark}[theorem]{Remark}
\newtheorem{remarks}[theorem]{Remarks}
\newtheorem{finalremarks}[theorem]{Final Remarks}
\newtheorem{definition}[theorem]{Definition}
\newtheorem{corollary}[theorem]{Corollary}
\newtheorem{example}[theorem]{Example}
\newtheorem{examples}[theorem]{Examples}
\newtheorem{assumption}[theorem]{Assumption}
\newtheorem{problem}[theorem]{Problem}
\newtheorem{question}[theorem]{Question}
\newtheorem{conjecture}[theorem]{Conjecture}
\newtheorem{rigiditytheorem}[theorem]{Rigidity Theorem}
\newtheorem{rigiditytheoremlevelN}[theorem]{Rigidity Theorem for the elliptic genus of level ${\boldsymbol {N}}$}
\newtheorem{maintheorem}[theorem]{Main Theorem}
\newtheorem{mainlemma}[theorem]{Main Lemma}
\newtheorem{fp criterion}[theorem]{Fixed point criterion}
\newtheorem{claim}[theorem]{Claim}
\newtheorem{illu}[theorem]{Proposition}

\newtheorem{strongrigiditytheorem}[theorem]{Strong rigidity theorem}
\newtheorem{ahatvanishingtheorem}[theorem]{${\boldsymbol {\hat A}}$-vanishing theorem}
\newtheorem{loopvanishingtheorem}[theorem]{${\boldsymbol {\hat A}}$-vanishing theorem for the loop space}
\newtheorem{bettinumbertheorem}[theorem]{Betti number theorem}
\newtheorem{stolzconjecture}[theorem]{Stolz' conjecture}
\newtheorem{frankelconjecture}[theorem]{Generalized Frankel conjecture}

\newtheorem{theoremsignrigid}[theorem]{Rigidity theorem for the signature}
\newtheorem{theoremToddrigid}[theorem]{Rigidity theorem for the Todd genus}
\newtheorem{theoremsignvanish}[theorem]{Vanishing theorem for the signature}
\newtheorem{theoremAdachvanish}[theorem]{Vanishing theorem for the ${\boldsymbol {\hat A}}$-genus}

\newtheorem{localization}[theorem]{Localization theorem for cohomology}
\newtheorem{localizationK}[theorem]{Localization theorem for ${\boldsymbol {K}}$-theory}
\newtheorem{fpf}[theorem]{Fixed point formula for cohomology}
\newtheorem{fpfK}[theorem]{Fixed point formula for ${\boldsymbol {K}}$-theory}
\newtheorem{fpfindK}[theorem]{Lefschetz fixed point formula}

\newcommand{\sslash}{\mathbin{/\mkern-6mu/}}

\def\Z{{\mathbb Z}}
\def\R{{\mathbb R}}
\def\Q{{\mathbb Q}}
\def\C{{\mathbb C}}
\def\N{{\mathbb N}}
\def\H{{\mathbb H}}
\def\Zp #1{{\mathbb Z }/#1{\mathbb Z}}

\def\sec{{\rm sec}}
\def\diam{{\rm diam}}

\def\ind{\mathrm {ind}\;}
\def\id{\mathrm {id}}
\def\im{\mathrm {im}}

\def\paperref#1#2#3#4#5#6{\text{#1:} #2, {\em #3} {\bf#4} (#5)#6}
\def\bookref#1#2#3#4#5#6{\text{#1:} {\em #2}, #3 #4 #5#6}
\def\preprintref#1#2#3#4{\text{#1:} #2 #3 (#4)}

\def\phiell#1{sign({\mathcal L}#1)}
\def\phinullell#1{\Phi_0(#1)}

\def\codim{{\rm{codim}\ }}

\def\sign{{\rm{sign}\ }}

\def\coker{{\rm {coker}\ }}

\def\dq{\slash\!\!\!\!\slash}

%%%%%%%%%%%%%%%%%%%%%%%%%%%%%%%%%%%%%%%%%%%%%%%%%%%%%%%%%%%%%%

\title[Torus actions, elliptic genera and positive curvature]{Torus Actions, Fixed-Point Formulas, Elliptic Genera and positive curvature}
\author{Anand Dessai}
\address{Department of Mathematics, Chemin du Mus\'ee 23, University of Fribourg, 1700 Fribourg (Switzerland) }
\email{anand.dessai@unifr.ch}
\urladdr{http://homeweb.unifr.ch/dessaia/pub/}
\date{}

\begin{abstract} We study fixed points of smooth torus actions on closed manifolds using fixed point formulas and equivariant elliptic genera. We also give applications to positively curved Riemannian manifolds with symmetry.
\end{abstract}

\maketitle

\noindent
\section{Introduction}\label{intro}
The purpose of this paper is to give a survey on some techniques used to study fixed points of smooth torus actions on closed manifolds. We will discuss and apply fixed point formulas as well as rigidity and vanishing theorems for classical operators and elliptic genera. These techniques have numerous applications in the study of differentiable transformation groups and, more recently, in the study of positively curved manifolds with symmetry.

The theory of transformation groups, even if one restricts to torus actions, is a vast field and our exposition is by no means complete. In particular, we will not focus on the question of existence of torus actions and will not describe techniques specific to manifolds with non-trivial fundamental group.

We hope that our survey will be useful for mathematicians working in the field and will fill a gap in the literature. For a discussion of some of the many other aspects we refer to \cite{Bo60, CoFl64, Br72, Ko72, Hs75, tD87, Ka91, AlPu93}.

In the first part of this paper we recall the localization theorems and fixed point formulas in cohomology, $K$-theory and index theory. The fixed point formulas are then applied in a more or less systematic way to describe various conditions under which a torus must act with fixed points. We also apply the fixed point formulas to give more precise information on the cohomology of the fixed point manifold for a large class of manifolds including cohomologically symplectic manifolds and manifolds with cohomology generated by classes of degree two. These results should be well-known to the experts as they only use techniques from the 1950s and 1960s. In the cases where we couldn't find a reference proofs are indicated for the convenience of the reader (see for example Theorem \ref{H2-orientable theorem}).

In the second part of this paper we discuss rigidity and vanishing theorems for classical operators and elliptic genera as well as applications to positively curved manifolds. We first recall the rigidity of classical operators, namely of the signature, Dirac and Dolbeault operator, and then describe the rigidity theorems for elliptic genera of level $2$ and higher level due to Witten, Taubes, Bott, Hirzebruch and others.

The rigidity under compact connected Lie group actions often forces related invariants to vanish. For elliptic genera, these invariants are given by the coefficients in the expansion in certain cusps. Their vanishing is related to finite cyclic subactions with fixed point manifold of large codimension and extends the classical vanishing in the cusps to vanishing of higher order. For example for the elliptic genus of level $2$ (the signature on the free loop space) of a spin manifold one obtains higher order versions of the $\hat A$-vanishing theorem of Atiyah and Hirzebruch. For (stably) almost complex manifolds with first Chern class divisible by $N\geq 2$ similar results hold for the elliptic genus of level $N$. This gives new information on the fixed point manifold. Another interesting consequence of the higher order vanishing is that it provides obstructions to isometric actions on manifolds of positive curvature.

To a large part this paper is of expository nature and most of the results are well-known, at least to the experts.  The higher order vanishing results for the elliptic genus of level $N$ presented in this paper and applications thereof to manifolds of positive curvature are new and have not been published.

The paper is structured as follows. In Section \ref{fpfsection} we review the localization theorem and fixed point formula for ordinary cohomology, topological $K$-theory and index theory. In Section \ref{basicsection} we apply the fixed point formulas to derive various fixed point theorems. Next to existence results for fixed points we will also discuss some structure results for the fixed point manifold of smooth torus actions. In Section \ref{classicaloperatorsection} we review the rigidity and vanishing theorems for the indices of classical operator and relations to positive/nonnegative curvature. In Sections \ref{level2section} and \ref{levelNsection} we recall the rigidity of elliptic genera of level $\geq 2$ and explain related higher vanishing results. These are then used to describe obstructions to isometric actions on manifolds of positive curvature.

\bigskip
\noindent

{\em Acknowledgements.} This survey article is an extended version of my talk at the conference {\em First Chinese-German Workshop on Metric Riemannian Geometry}. It is my pleasure to thank the organizers and the participants for the interesting and very enjoyable conference. I would also like to thank the referee for helpful remarks.

\section{Localization and fixed point formulas}\label{fpfsection}
In this section we will review the localization theorems and fixed point formulas for ordinary cohomology, topological $K$-theory and index theory. These formulas will be applied in the next sections to derive various fixed point theorems, to study the structure of fixed point sets, to discuss rigidity and vanishing theorems for elliptic genera and to describe applications for manifolds of positive curvature.

We begin by recalling the Borel construction. Let $G$ be a compact Lie group, $EG\to BG$ a classifying bundle for $G$ and $X$ a $G$-space. Note that $G$ acts freely on $EG\times X$ via the diagonal action and $EG\times X$ is homotopically equivalent to $X$. Following Borel  \cite{Bo60}, one associates to the $G$-space $X$ the fiber bundle \linebreak  $\pi :X_G\to BG$ with fiber $X$ where $X_G:=EG\times _G X:=(EG\times X)/G$ is the quotient by the free $G$-action. The space $X_G$ is called the {\em Borel construction} of $X$ or {\em homotopy quotient} of $X$ by $G$.

The use of this construction is motivated by the idea that the behavior of the $G$-action is partially reflected in the cohomological properties of the fiber bundle \linebreak $\pi :X_G\to BG$ \wrt \ a chosen cohomology theory, e.g. ordinary cohomology, topological $K$-theory, cobordism etc. For example the existence of a fixed point $pt$ of the $G$-action implies that the homomorphism in cohomology induced by $\pi $ is injective since the equivariant map $j:pt\to X$ induces a section $j_G:pt_G=BG\to X_G$ of the projection $\pi$. This approach to group actions has been developed in Borel's seminar on transformation groups in the late 1950s \cite{Bo60}. For de Rham cohomology a different construction using differential forms was already carried out by Cartan in the early 1950s (see \cite{GuSt99}).

An important relation between the cohomological properties of the fiber bundle $X_G\to BG$ and the fixed point set is expressed by the localization theorem. Roughly speaking this theorem says that after localization the inclusion of the fixed point set $X^G\hookrightarrow X$ induces an isomorphism in the equivariant cohomology theory under consideration. The first example of this kind was Borel's localization theorem for cohomology \cite[Ch. IV]{Bo60}. A few years later, after topological $K$-theory was introduced by Atiyah and Hirzebruch \cite{AtHi59,AtHi61}, the localization theorem for $K$-theory was established by Atiyah and Segal \cite{AtSeII68,Se68}.

Localization theorems are known to hold for a wide range of equivariant multiplicative cohomology theories (see for example \cite{tD70,Qu71a,Qu71b}, \cite[Ch. III.3]{tD87}, \cite[Ch. 6.2]{Ka91}). Here we will restrict to ordinary cohomology and topolo\-gical $K$-theory. Also we will formulate the localization theorems only in the case of smooth torus actions. This will be sufficient for the purposes of this paper.

The localization theorems lead to corresponding fixed point formulas. In equivariant $K$-theory and index theory these formulas were established in the 1960s by Atiyah, Bott, Segal and Singer and were used at that time to derive various fixed point formulas in cohomology for equivariant characteristic classes \cite{AtBoII68, AtSiIII68}. Bott's residue formula provides another approach to cohomological fixed point formulas via Chern-Weil theory \cite{Bo67a, Bo67b} (see also \cite[Ch. II, \S 6]{Ko72}). The full fixed point formulas for ordinary cohomology apparently were not considered prior to the influential work of Atiyah-Bott and Berline-Vergne on their de Rham version of the fixed point formula in the early 1980s \cite{AtBo84,BeVe82}.

\subsection{Cohomology}
Let $T$ be a torus of rank $r$ and let $S$ be the multiplicative set of non-zero cohomology classes in $H^*(BT;\Z)$. For any $H^{*}(BT;\Z)$-module $A$ let $S^{-1}A$ be its localization \wrt \ $S$. Note that if we identify $H^{*}(BT;\Z)$ with $\Z [u_1,\ldots ,u_r]$, $u_1,\ldots ,u_r\in H^2(BT;\Z)$ a basis, then $S^{-1}\Z [u_1,\ldots ,u_r]$ is the quotient field $\Q (u_1,\ldots ,u_r)$ and $S^{-1}A$ is identified with $A\otimes _{\Z [u_1,\ldots ,u_r]}\Q (u_1,\ldots ,u_r)$.

For any $T$-space $X$ let $H^*_{T}(X ;\Z ):=H^*(X_{T};\Z)$. We apply corresponding notations for other coefficients and other group actions. Note that $H^*_{T}(X ;\Z )$ is a $H^{*}(BT;\Z)$-module \wrt \ the module action induced by the projection $X_T\to BT$. We can now state the localization theorem (see \cite[Ch. IV]{Bo60}, \cite[Thm. 4.4]{Qu71a}, \cite[Ch. VII]{Br72}, \cite[Ch. III, \S 2]{Hs75}, \cite[Ch. 3]{AlPu93}).

\begin{localization} Let $M$ be a manifold with smooth $T$-action. Then the inclusion of the fixed point manifold $M^{T}$ induces an isomorphism $$S^{-1}H^*_{T}(M;\Z )\overset \cong \to  S^{-1}H^*_{T}(M^{T};\Z ).$$\proofend
\end{localization}

Note that since $M^{T}$ is a trivial $T$-space we have $(M^{T})_{T}= M^{T}\times BT$ and $H^*_{T}(M^{T};\Z )\cong H^*(M^{T};\Z )\otimes H^*(BT;\Z )$ by the K\"unneth theorem.

The localization theorem above basically says that the torsion-free part of the $H^{*}(BT;\Q)$-module $H^*_{T}(M;\Q )$ is geometrically carried by the fixed point manifold $M^{T}$. Also, if $M$ is $n$-dimensional then $H^i_{T}(M;\Q )\to H^i_{T}(M^{T};\Q )$ is an isomorphism for $i>n$ (see for example \cite[Ch. VII, Thm. 1.5]{Br72}, for a refined version of the localization theorem see \cite[Thm. 4.4]{Qu71a}).

The proof of the localization theorem and the properties mentioned above uses the Leray-Serre spectral sequence for $M_{T}\to BT$ together with the fact that $H^*_{T}(M-M^{T};\Q )$ is a torsion $H^{*}(BT;\Q)$-module. Corresponding results hold for elementary $p$-abelian groups (see \cite[Ch. IV]{Bo60}, \cite[Ch. VII]{Br72}, \cite[Ch. III, \S 2]{Hs75}, \cite[Ch. 3]{AlPu93}). If one describes $S^{-1}H^*_{T}(M;\Z )$ in terms of generators and relations then the number of connected components and their localized cohomology can be computed in principle from the ideal of defining relations (see \cite[Ch. IV, \S 1]{Hs75} for a precise statement).

From the localization theorem one readily sees the following fixed point criterion:
$$M^{T}\neq\emptyset \iff  H^*_{T}(M;\Q ) \text{ is not a torsion } H^*(BT;\Q )\text{-module}.$$

We now turn to the fixed point formula for ordinary cohomology. Let us first recall that ordinary integral cohomology is oriented \wrt \ oriented vector bundles. Hence, an oriented vector bundle comes with a Thom class and a Thom isomorphism \cite{Th52,MiSt74}. This allows to define a push-forward $f_!$ (also called Gysin homomorphism or umkehr homomorphism) for an oriented map $f:X\to Y$ between closed smooth manifolds (see for example \cite{AtBo84}).\footnote{The map $f$ is called {\em oriented} if the bundle $TX\ominus f^*(TY)$ is oriented (for the general concept of orientability see for example \cite{Sw75}, Ch. 14]).} Let $k$ denote the difference of the dimensions of $Y$ and $X$. Then $f_!$ is a $H^{*}(Y;\Z )$-module homomorphism from $H^*(X;\Z )$ to $H^{*+k}(Y;\Z )$ where $H^{*}(Y;\Z )$ acts on $H^*(X;\Z )$ via $f^*$. If $X$ and $Y$ are oriented then $f_!$ is given by pre- and post-composing the induced homomorphism $f_*:H_*(X;\Z )\to H_*(Y;\Z )$ in homology with the Poincar\'e isomorphisms for $X$ and $Y$.

Similarly one can consider push-forward maps in the equivariant setting using the Borel construction (see for example \cite{AtBo84}). We will not give the details on the construction here and will only mention the following two basic cases for an $n$-dimensional oriented closed manifold $M$ with a smooth action by the torus $T$.

In the first case we consider the projection of $M$ to a point $pt$. The projection induces a push-forward
$$\int _M:H^*_{T}(M;\Z )\to H^{*-n}_{T}(pt;\Z )\cong H^{*-n}(BT;\Z ).$$
If the $T$-action is trivial then $\int _M$ is given by evaluation on the fundamental cycle of $M$.

In the second case we consider a connected component $Y$ of codimension $k$ of the fixed point manifold $M^T$. Let $j_Y:Y\hookrightarrow M$ denote the inclusion and let $\nu _Y$ be the $T$-equivariant normal bundle of $Y$. Note that $Y$ and $\nu _Y$ are orientable and orientations for them determine each other if one requires compatibility with the given orientation of $M$. For fixed orientations the inclusion $j_Y$  induces a push-forward $(j_Y)_!:H^*_{T}(Y;\Z )\to H^{*+k}_{T}(M;\Z )$ which is given by multiplication with the equivariant Thom class of $\nu _Y$. For later reference we recall that the restriction of this class to $Y$ is the equivariant Euler class $e_{T}(\nu _Y)\in H^k_T(Y;\Z)$. We are now ready to state the

\begin{fpf}\label{fpf} Let $M$ be an $n$-dimensional oriented closed manifold with a smooth action by the torus $T$. Let $v\in H^*_T(M;\Z )$. Then
$$\int _M v=\sum _{Y\subset M^T} \int _Y \frac {j_Y^*(v)} {e_T(\nu _Y)}.$$
\end{fpf}

\bigskip
\noindent
Note that the left hand side is an element in $H^{*}(BT;\Z )$ whereas each summand on the right hand side is an element in the localized ring $S^{-1}H^{*}(BT;\Z )$. Since the homomorphism $H^{*}(BT;\Z )\to S^{-1}H^{*}(BT;\Z )$ is injective the identity in the fixed point formula above holds in $S^{-1}H^{*}(BT;\Z )$ and also in $H^{*}(BT;\Z )$. For the convenience of the reader we recall the proof.

\bigskip
\noindent
{\bf Proof:} The fixed point formula follows from properties of the push-forward and the localization formula. Let $j:M^{T}\hookrightarrow M$ be the inclusion of the fixed point manifold and let $j_!$ denote the push-forward for a fixed orientation of $M^{T}$. Similarly let $j_Y:Y\hookrightarrow M$ denote the inclusion of a fixed point component of codimension $k$ and let $(j_Y)_!$ be the push-forward \wrt \ the orientation of $Y\subset M^T$. By functoriality  $\int _{Y}$ factorizes as $\int _M\circ (j_Y)_!$. The composition
$$j_Y^*\circ (j_Y)_!:H_T^*(Y;\Z )\to H_T^{*+k}(Y;\Z )$$ is given by multiplication with $e_T(\nu _Y)$. Since $T$ acts without fixed points outside of the zero section of $\nu _Y$ the restriction of $e_T(\nu _Y)$ to a point of $Y$ is non-zero. This together with the nilpotency of  $H^{>0}(Y;\Z )$ implies that $e_T(\nu _Y)$ is invertible in the localization $S^{-1}H^*_{T}(Y;\Z )$, where $S$ denotes the multiplicative set of non-zero elements in $H^{*}(BT;\Z)$.

Recall from the localization theorem that $j^*:S^{-1}H^*_{T}(M;\Z ) \to  S^{-1}H^*_{T}(M^{T};\Z )$ is an isomorphism. Hence, $$v=\sum_{Y\subset M^T} (j_Y)_!\left (\frac {j_Y^*(v)}{e_T(\nu _Y)}\right )\in S^{-1}H^*_{T}(M;\Z )$$
 since both sides have the same image under $j^*$. This gives
$$\int _M v=\sum_{Y\subset M^T} \int_M (j_Y)_!\left (\frac {j_Y^*(v)}{e_T(\nu _Y)}\right )=\sum _{Y\subset M^T} \int _Y \frac {j_Y^*(v)} {e_T(\nu _Y)}.$$\proofend

In the case of a circle action with isolated fixed points the fixed point formula for equivariant characteristic numbers takes the form of Bott's residue formula \cite{Bo67a}. Other cases of this formula were formulated at that time using the fixed point formulas for equivariant $K$-theory and equivariant index theory \cite{AtBoII68, AtSiIII68}. A de Rham version of Fixed point formula \ref{fpf} is due to Atiyah and Bott \cite{AtBo84} as well as Berline and Vergne \cite{BeVe82} (see also \cite{Tu11} for the historical development).

\subsection{${\boldsymbol {K}}$-theory and index theory}
As mentioned before localization theorems exist not only for ordinary cohomology but can be established for a wide class of multiplicative cohomology theories. Of fundamental importance are the localization theorem for topological $K$-theory and the fixed point formulas for $K$-theory and index theory which we will briefly recall.

The $K$-theory $K(X)$ of a compact Hausdorff space $X$ is defined by turning the semi-ring of isomorphism classes of complex vector bundles over $X$ into a ring via the Grothendieck construction \cite{AtHi59,AtHi61,At67}. If $X$ is a locally compact Hausdorff space then  $K(X)$ shall denote $K$-theory with compact support (see \cite{AtSeII68,Se68}).\footnote{The $K$-theory with compact support of a locally compact Hausdorff space is isomorphic to the reduced $K$-theory of its one-point compactification.}

This construction refines directly in the presence of a group action. In fact, if a compact Lie group $G$ acts on $X$ then - by considering $G$-equivariant complex vector bundles over $X$ - one obtains the equivariant $K$-theory of $X$, $K_G(X)$, which is a finer invariant than $K(X_G)$ by the Atiyah-Segal completion theorem. For example, for a point $pt$ the equivariant $K$-theory $K_G(pt)$ is isomorphic to the complex representation ring $R(G)$ whereas $K(pt_G)\cong K(BG)$ is the completion of $R(G)$ \wrt \ the augmentation ideal \cite{AtSe69}.

The localization theorem for $K$-theory takes for an action by a torus $T$ the following form (see \cite{Se68}, \cite[Localization theorem 1.1]{AtSeII68}).
\begin{localizationK} Let $X$ be a locally compact $T$-space and let $i:X^{T}\hookrightarrow X$ denote the inclusion of the fixed point set. Then the homomorphism induced by the inclusion
 $$i^*: K_{T}(X)\to  K_{T}(X^{T})$$ becomes an isomorphism
when localized \wrt \ the multiplicative set of non-zero elements of $R(T)$.\proofend
\end{localizationK}

We now turn to the fixed point formula for $K$-theory. The reasoning is analogous to the one for cohomology. We will skip the details since we will mainly use in the following sections the index theoretical version \ref{fpfindK} (see  \cite[Lemma 2.5]{AtSeII68} for details).

Let us first recall that $K$-theory is oriented for complex vector bundles (or, more generally, for spin$^c$ vector bundles) \cite{AtBoSh64}. This property allows to define push-forward maps. We will only need the following two basic cases for an almost complex manifold $M$ with a smooth action by a torus $T$ preserving the structure.

In the first case consider the projection $\pi ^M:M\to pt$. The projection induces a push-forward
$$\pi^M_!:K_{T}(M)\to K_T(pt)\cong R(T)$$ which is a $K_T(pt)$-module homomorphism.

In the second case we consider a connected component $Y$ of the fixed point manifold $M^T$. Let $\nu _Y$ be the $T$-equivariant normal bundle of $Y$. Note that the almost complex structure on $M$ induces an almost complex structure on $Y$ and a $T$-equivariant complex structure on $\nu _Y$.

With respect to the complex structure of the normal bundle the inclusion \linebreak $j_Y:Y\hookrightarrow M$ induces a push-forward $(j_Y)_!:K_{T}(Y)\to K_{T}(M)$ which is a $K_{T}(M)$-module homomorphism \wrt \ the action of $K_{T}(M)$ on $K_{T}(Y)$ via $j_Y^*$.

The push-forward $(j_Y)_!$ is given by multiplication with the equivariant $K$-theore\-tical Thom class of $\nu _Y$. For later reference we recall that the restriction of this class to $Y$ is the equivariant $K$-theoretical Euler class $\Lambda _{-1}(\nu _Y)\in K_T(Y)$, where $\Lambda _t(\nu _Y):=\sum _{i\geq 0} \Lambda ^i(\nu _Y) t^i$ and $\Lambda^i(\nu _Y)$ is the $i$th exterior power of $\nu _Y$. We are now ready to state the fixed point formula for $v\in K_T(M)$:
$$\pi^M_! (v)=\sum _{Y\subset M^T} \pi^Y_! \left (\frac {j_Y^*(v)} {\Lambda _{-1}(\nu _Y)}\right )$$
The formula also holds under the weaker condition that $M$ has a $T$-equivariant spin$^c$ structure. The $K$-theoretical Euler class $\Lambda _{-1}(\nu _Y)$ is non-zero in $R(T)$ after restriction to any point of $Y$. This implies that $\Lambda _{-1}(\nu _Y)$ is invertible in the localization (see \cite[Lemma 2.5]{AtSeII68}). The proof now follows the same line of reasoning as the proof of Fixed point formula \ref{fpf}.

Next we will recall the fixed point formula for equivariant index theory. This will be our main tool when we discuss rigidity and vanishing theorems for elliptic genera and applications thereof to group actions.

Let $M$ be a closed $n$-dimensional smooth manifold with tangent bundle\linebreak $\pi ^M:TM\to M$. Note that the tangent bundle of $TM$ is isomorphic to the pull back via $\pi ^M$ of $TM\oplus TM$ and the manifold $TM$ comes with an almost complex structure. This allows to define the push-forward map \wrt \ the map $\pi ^{TM}:=(\pi ^M)_*:TM\to Tpt\cong pt$:
$$\pi ^{TM}_!:K(TM)\to K(pt)\cong \Z.$$

Now let $D$ be an elliptic differential operator on $M$, resp. an elliptic complex
$$\ldots \overset {D_{i-1}}\longrightarrow \Gamma (E_i) \overset {D_{i}}\longrightarrow \Gamma (E_{i+1})\overset {D_{i+1}}\longrightarrow \Gamma (E_{i+2})\to \ldots ,$$
with symbol $\sigma (D)\in K(TM)$ and index $\ind (D)=\dim (\ker (D)) - \dim (\coker (D))$, resp. $\ind (D)=\sum _i (-1)^i \dim (\ker (D_i)/\im (D_{i-1}))\in \Z$. By the celebrated Atiyah-Singer index theorem \cite{AtSiI68} the index of $D$ can be computed by applying the push-forward to the symbol of $D$, i.e.
$$\ind (D)=\pi ^{TM}_!(\sigma(D)).$$

Prominent applications of the index theorem are to the signature operator for orientable manifolds, the Dolbeault operator for almost complex manifolds and the Dirac operator for spin manifolds which lead to the signature theorem, the Riemann-Roch theorem and the integrality theorem for the $\hat A$-genus, respectively (see \cite{AtSiIII68} for details).

Next we recall the equivariant index theorem. Suppose a compact Lie group $G$ acts smoothly on $M$ and $D$ is a $G$-equivariant elliptic differential operator, resp. $G$-equivariant elliptic complex. In this situation the index $\ind (D)$ is an element of the representation ring $R(G)$, the symbol $\sigma (D)$ is an element of the equivariant $K$-theory $K_G(TM)$ and the index theorem \cite{AtSiI68} asserts that
$$\ind (D)=\pi ^{TM}_!(\sigma(D))\in R(G),$$
where $\pi ^{TM}_!:K_G(TM)\to K_G(pt)\cong R(G)$ is the push-forward for the $G$-equivariant map $\pi ^{TM}:TM\to Tpt\cong pt$.

We now come to the fixed point formula for index theory. Let $G$ be topologically cyclic generated by $g\in G$. Let $j_Y:Y\hookrightarrow M$ be the inclusion of a fixed point component $Y\subset M^g=M^G$ and let $\nu_Y$ be the normal bundle. We note that the normal bundle of $j:TY\hookrightarrow TM$ comes with a complex structure which allows to define a push-forward $K_G(TY)\to K_G(TM)$. The normal bundle together with this complex structure defines the element $\nu_Y\otimes \C\in K_G(TY)$. We are now ready to state the fixed point formula for equivariant index theory (see \cite[Prop. 2.8]{AtSeII68}, \cite[Prop. 2.10, Thm. 2.12]{AtSe69}).

\begin{fpfindK}\label{fpfindK} Let $M$ be a closed manifold with smooth $G$-action, were $G$ is topologically cyclic generated by $g\in G$. Let $D$ be a $G$-equivariant elliptic differential operator (or complex) on $M$ with symbol $u:=\sigma (D)$. Then its index $\ind (D)\in R(G)$ evaluated at $g$ is given by
$$\ind (D)(g)=\sum _{Y\subset M^g} a_Y,$$
where
$$a_Y=\pi ^{TY}_!\left (\frac {u(g)_{\vert TY}}{\Lambda _{-1}(\nu_Y\otimes \C ) (g)}\right ).$$\proofend
\end{fpfindK}

Here $u(g)_{\vert TY}$ denotes the restriction of $u(g)$ to $TY$, i.e. the image under the restriction homomorphism $K(TM)\otimes \C \to K(TY)\otimes \C $. The Lefschetz fixed point formula can also be phrased in terms of rational cohomology using the natural transformation from $K$-theory to cohomology via the Chern character (see \cite[Lefschetz Theorem 3.9]{AtSiIII68}). The cohomological version, which will be described in Subsection \ref{index theorem in cohomology}, will be our mail tool when we discuss rigidity and vanishing theorems for elliptic genera.

\section{Fixed points}\label{basicsection}
In this section we will apply Fixed point formula \ref{fpf} to derive various fixed point theorems. All the results presented here can be equally well shown using the fixed point formulas for $K$-theory or index theory. Next to existence results for fixed points we will also discuss some structure results for the fixed point manifold. As before we will restrict to smooth torus actions. From now on all manifolds will be assumed to be connected.

\subsection{Existence of fixed points}
Let $M$ be a closed oriented $n$-dimensional smooth manifold and let $T$ be a torus of rank $r$ which acts smoothly on $M$. Suppose the fixed point set $M^T$ is empty. Then, by Fixed point formula \ref{fpf} for any class $v\in H^*_T(M;\R )$ one has $\int _M v=0$. For further reference we single out this classical observation in the following form:

\begin{fp criterion}\label{fp criterion} Suppose there exists an equivariant class $v\in H^i_T(M;\R )$ such that $\int _M v \in H^{i-n}_T(pt;\R)=H^{i-n}(BT;\R)$ is non-zero. Then $M^T$ is non-empty.\end{fp criterion}

\bigskip
\noindent

The following well-known special case already accounts for many applications.

\begin{theorem}\label{basicfpt}
Suppose there exists a cohomology class $w$ in the image of the restriction map $H^n_T(M;\R )\to H^n(M;\R )$ such that $\int _M w \neq 0\in \R$. Then $M^T$ is non-empty.
\end{theorem}

\pproof Let $v$ be a pre-image of $w$ \wrt \ $H^n_T(M;\R )\to H^n(M;\R )$. Since the restriction map commutes with integration and $H^0_T(pt;\R )\to H^0(pt;\R )$ is an isomorphism it follows that $\int _M v=\int _Mw\neq 0$. Hence, $M^T$ is non-empty.\proofend

Let us say $w\in H^*(M;\R )$ {\em has an equivariant extension} if there exists a class $v\in H^*_T(M;\R )$ which maps to $w$ under the restriction map $H^*_T(M;\R )\to H^*(M;\R )$.

Since the $T$-action on $M$ extends to an action on the tangent bundle $TM$ via differentials all characteristic classes of $M$ have equivariant extensions. Hence, one can apply the last theorem to the Euler class (unstable characteristic class) and to the Pontrjagin classes (stable characteristic classes) of $M$.

Taking $w$ to be the Euler class $e(M)$ and noting that $\int _M e(M)$ is the Euler characteristic $\chi (M)$ one obtains the following fixed point theorem which is a special case of a classical result of Hopf on the zeros of vector fields \cite{Ho26}.

\begin{illu} Let $M$ be an oriented closed manifold with $\chi (M)\neq 0$. Then any smooth action by a torus $T$ on $M$ has at least one fixed point.\proofend
\end{illu}

The orientability condition here is not necessary. If $M$ is not orientable one can lift a two-fold action to the orientation cover and conclude that the $T$-action on $M$ has a fixed point.

Taking $w\in H^n(M;\Z )$ to be a polynomial in the Pontrjagin classes of $M$ so that $\int _M w$ is a Pontrjagin number one obtains the following classical result (see for example \cite[Ch. II, \S 6]{Ko72}).

\begin{illu} Let $M$ be an oriented closed manifold for which not all Pontrjagin numbers vanish. Then any smooth action by a torus $T$ on $M$ has at least one fixed point.
\proofend
\end{illu}

As shown by Thom [Th54, MiSt74] all torsion in the oriented bordism ring is finite and a closed oriented manifold represents a torsion class in the oriented bordism ring if and only if all its Pontrjagin numbers vanish. Hence, the last proposition implies the following well-known

\begin{corollary} If $S^1$ acts smoothly on $M$ without fixed points then $M$ is rationally zero in the oriented bordism ring, i.e. a non-zero multiple of $M$ is the boundary of an oriented $(n+1)$-dimensional compact manifold.\proofend\end{corollary}

The result above can be refined in the presence of a reduction of the stable structure group of $M$. To illustrate this let us take a look at stably almost complex manifolds. Recall that a stably almost complex structure on $M$ is given by a complex structure of its stable tangent bundle, i.e. an endomorphism $J$ of $TM\oplus \epsilon_k$ with $J^2=-\id$ (here $\epsilon_k$ denotes the trivial $k$-dimensional real vector bundle over $M$). The manifold $M$ together with a stably almost complex structure is called a stably almost complex manifold or unitary manifold. Note that such a structure defines an orientation on $M$ and (after fixing a $J$-equivariant Riemannian metric) a reduction of the stable structure group to $U((n+k)/2)\hookrightarrow SO(n+k)$. Now applying Theorem \ref{basicfpt} to the equivariant Chern numbers gives the following classical result.

\begin{illu}\label{Chernumber fixedpoint} Let $M$ be a stably almost complex manifold for which not all Chern numbers vanish. Then any smooth action by a torus $T$ on $M$ which lifts to the stably almost complex structure has at least one fixed point.\proofend
\end{illu}

From the work of Milnor \cite{Mi60} and Novikov \cite{No60} one knows that a stably almost complex manifold represents the zero element in the unitary bordism ring if and only if all its Chern numbers vanish. Hence, the last proposition implies the well-known

\begin{corollary} Let $M$ be a stably almost complex manifold. Suppose $S^1$ acts smoothly on $M$ and the action lifts to the stably almost complex structure. If $S^1$ acts without fixed points then $M$ is zero in the complex bordism ring, i.e. $M$ is the boundary of an $(n+1)$-dimensional stably almost complex compact manifold.\proofend\end{corollary}

The situation is quite similar for other reductions of the structure group. Here we only point out the following fixed point theorem. Let $H$ be a compact Lie group. Suppose $M$ admits a stable reduction of its structure group to $H$. More precisely, we have an $H$-principal bundle $P\to M$ together with a homomorphism of Lie groups $\rho: H\to SO(n+k)$  and an isomorphism of oriented vector bundles $P\times_\rho \R ^{n+k}\overset \cong \longrightarrow TM\oplus \epsilon_k$. Now let $f:M\to BH$ be a classifying map for $P\to M$ and consider the (real) characteristic ring $f^*(H^*(BH;\R ))\subset H^*(M;\R )$. The {\em characteristic numbers} of $P$ are defined by $\int _M v$ for $v\in f^*(H^n(BH;\R ))$.

Let $T$ act smoothly on $M$ and linearly on $\epsilon_k$. Now suppose the action lifts to $P$, i.e. $P\to M$ is a $T$-equivariant principal bundle and the isomorphism $P\times_\rho \R ^{n+k}\overset \cong \longrightarrow TM\oplus \epsilon_k$ is equivariant \wrt \ the $T$-action. Then the characteristic ring is in the image of $H^*_T(M;\R )\to H^*(M;\R )$. Arguing as before one gets the following generalization of Proposition \ref{Chernumber fixedpoint}.

\begin{illu}\label{reduction structure group fixed point} Let $M$ be an oriented closed manifold and $P\to M$ a stable reduction of its structure group for which not all characteristic numbers of $P$ vanish. Then any smooth $T$-action on $M$ which lifts to $P$ has at least one fixed point.\proofend
\end{illu}

We remark that this result can be generalized further from equivariant stable reductions to equivariant $B$-structures (see \cite{La63} for the definition of a $B$-structure). We will not give the details here. A particular case of such a fixed point theorem is given in the next theorem.

Note however that it is in general difficult to decide whether the $T$-action on $M$ lifts to the principal bundle $P\to M$. Also the question whether a cohomology class $w\in H^*(M;\R )$ has an equivariant extension, i.e. is in the image of the homomorphism $H^*_T(M;\R )\to H^*(M;\R )$, has often to be left open. A simple well-known criterion for a positive answer is given in the next lemma which follows from an inspection of the Leray-Serre spectral sequence for the Borel-construction $M\hookrightarrow M_T\to BT$.
\begin{lemma} If $b_{2i+1}(M)=0$ for $0<2i+1<2k$ then any cohomology class of degree $2k$ has an equivariant lift.\proofend\end{lemma}

In particular, any cohomology class of degree $2$ has an equivariant extension if $b_1(M)=0$, e.g. if $M$ is simply connected. This already leads to some interesting fixed point theorems. For oriented manifolds, for example, one has:

\begin{proposition} Let $M$ be an oriented closed manifold with $b_1(M)=0$. Suppose there exist cohomology classes $x_1,\ldots ,x_l \in H^2(M;\Z)$ and Pontrjagin classes $p_{i_1}(M),\ldots , p_{i_s}(M)$ such that
$\int _M \prod_{i=1}^l x_i\cdot \prod _{j=1}^s p_{i_j}(M)\neq 0$. Then any smooth action by a torus $T$ on $M$ has at least one fixed point.
\end{proposition}

\pproof Let $w:=\prod_{i=1}^l x_i\cdot \prod _{j=1}^s p_{i_j}(M)$. Since $b_1(M)=0$ the classes $x_i$ admit equivariant extensions. It follows that the same is true for $w$. Since $\int _M w\neq 0$ the fixed point manifold $M^T$ is not empty by Theorem \ref{basicfpt}.\proofend

If $b_1(M)\neq 0$ then the conclusion holds if there exist equivariant extensions of the classes $x_i$. We invite the interested reader to formulate corresponding fixed point theorems in the presence of an equivariant reduction of the structure group and/or if $b_{2i+1}(M)=0$ for $0<2i+1<2k$ and $k>0$.

\subsection{Structure of ${\boldsymbol {M^T}}$}
In the remaining part of this section we will apply the cohomological fixed point formula to give some information on the structure of the fixed point manifold for a large class of manifolds including cohomologically symplectic manifolds and manifolds with cohomology generated by classes of degree two.

A general theorem which relates the localized equivariant cohomology of the manifold to the localized equivariant cohomology of the individual fixed point components is due to Hsiang \cite[Thm. IV.1, p. 47]{Hs75}. This theorem can be nicely applied to recover the structure results for the fixed point set of cohomology spheres or cohomology projective spaces (see \cite[Ch. IV, \S 1]{Hs75}). For more complicated spaces other techniques like the fixed point formulas are often easier to use to obtain information on the structure of the fixed point components. This will be illustrated in Theorem \ref{H2-orientable theorem}. As before we will restrict to smooth actions to simplify the exposition and refer to \cite{Br72,Hs75} for generalizations.

As a partial motivation for the following let us first recall the classical structure theorem for circle actions on cohomology complex projective spaces.

Let $M$ be an oriented closed $2m$-dimensional manifold with smooth $S^1$-action. Suppose $M$ is an integral cohomology complex projective space (i.e. the integral cohomology ring of $M$ is isomorphic to the integral cohomology ring of the complex projective space $\C P^m$). Let $x\in H^2(M;\Z )$ be a generator, let $Y_1,\ldots ,Y_l\subset M^{S^1}$ be the connected fixed point components and let $2m_i:=\dim Y_i$.

Then one has the following information for the fixed point components: Each $Y_i$ is an integral cohomology $\C P ^{m_i}$, the restriction of $x$ to $Y_i$ is a generator of $H^2(Y_i;\Z )$ and $\sum _i (m_i +1) =m+1$. If $M$ is a rational or real cohomology complex projective space then the statement holds after replacing integer coefficients by rational or real coefficients, respectively (see \cite[Ch. VII, \S 2]{Br72}, \cite[Thm. IV.3]{Hs75}). For a $K$-theoretical approach see \cite[Ch. II]{Pe72}.

As we will explain certain aspects of this structure result can be carried over to many other manifolds using Fixed point formula \ref{fpf}.

For the applications in the following sections it will be convenient to rather work with complex line bundles than with cohomology classes of degree two. To begin with we will briefly discuss lifting of group actions to complex line bundles.

Let $M$ be an oriented closed manifold on which a torus $T$ acts smoothly. Let $L\to M$ be a complex line bundle. Suppose the $T$-action lifts to $L$. For a fixed lift the restriction of the equivariant complex line bundle $L$ to a point of a fixed point component $Y\subset M^T$ reduces to a complex one dimensional $T$-representation. The $T$-weight $w_Y$ of this representation is independent of the point in $Y$ but may depend on $Y$. We call $w_Y$ the {\em local weight} of $L$ at $Y$. We may change the torus action on $L$ by taking the tensor product $L\otimes \chi $, where $\chi $ is a one-dimensional complex $T$-representation of weight $w$. Note that $L\otimes \chi $ and $L$ are canonically isomorphic as non-equivariant line bundles. For the new $T$-action the local weight at $Y$ is $w_Y+w$. Hence, the local weights of $L$ change simultaneously by a global weight $w$. In particular, we may choose for a given fixed point component $Y\subset M^T$ the lift of the $T$-action to $L$ such that the local weight at $Y$ vanishes, i.e. such that $T$ acts trivially on the restriction $L_{\vert Y}$. Conversely, any two lifts of the $T$-action to $L$ differ by a global weight.

The question whether an action lifts to $L$ has been studied by Stewart, Su, Hattori-Yoshida and others. It turns out that the only obstruction to a lift is the necessary condition that the first Chern class $c_1(L)$ of $L$ has an equivariant extension (see \cite{St61,Su63,HaYo76}). As pointed out before this is for example the case if $b_1(M)=0$.

We now come to an application of Fixed point formula \ref{fpf} which generalizes certain aspects of the structure result for cohomology complex projective spaces to many other manifolds.

Let us say that a closed connected $2m$-dimensional oriented manifold {\em $M$ can be oriented by degree $2$ classes} or {\em $M$ is $H^2$-orientable}, for short, if there are classes $x_1,\ldots ,x_m\in H^2(M;\Z )$ such that $\int _M (x_1\cdot \ldots \cdot x_m)\neq 0$, i.e. if $x_1\cdot \ldots \cdot x_m$ defines an orientation on $M$. We note that this class of manifolds includes for example cohomologically symplectic manifolds or oriented closed manifolds for which the real cohomology ring is generated by classes of degree $2$ (e.g. manifolds of the cohomological type of a quasitoric manifold).

In the situation where a torus $T$ acts on $M$ we say {\em $M$ is equivariantly $H^2$-orientable} if $M$ is $H^2$-orientable and the classes $x_1,\ldots ,x_m$ all have an equivariant extension.

\begin{theorem}\label{H2-orientable theorem} Let $M$ be an oriented closed connected $2m$-dimensional manifold. Let $T$ be a torus which acts smoothly on $M$ with fixed point components $Y_1,\ldots ,Y_l$ and $\dim Y_i=2m_i$. Suppose $M$ is equivariantly $H^2$-orientable. Then $m+1\leq \sum _{i=1}^l(m_i+1)$. Moreover, if equality holds then each component $Y_i$ is $H^2$-orientable.
\end{theorem}

This result should be well-known to the experts as the proof only uses techniques from the 1960s. Since we couldn't find a reference we give a proof for the convenience of the reader (see also \cite[\S 4]{DeWi}).

\bigskip
\noindent
\pproof Let $x_1,\ldots ,x_m\in H^2(M;\Z )$ be classes with $\int _M (x_1\cdot \ldots \cdot x_m)\neq 0$ such that each $x_i$, $i=1,\ldots ,m$, has an equivariant extension. Let $L_i$ be the complex line bundle over $M$ with $c_1(L_i)=x_i$, $i=1,\ldots ,m$. From the lifting properties discussed above we conclude that the $T$-action lifts to each of the line bundles.

Now assume that $m+1\geq \sum _{i=1}^l (m_i+1)$. We will choose equivariant extensions of the classes $x_1,\ldots ,x_m$ such that at least $l-1$ summands in the fixed point formula for $\int _M (x_1\cdot \ldots \cdot x_m)$ will vanish.

We adjust the lift of the $T$-action to the first $(m_1+1)$ line bundles such that $T$ acts trivially on $L_1, \ldots ,L_{m_1+1}$ after restriction to $Y_1$. In other words we choose the lift such that the local weights at $Y_1$ vanish for the first $(m_1+1)$ line bundles. Next we choose a lift to the $(m_2+1)$ line bundles $L_{m_1+2}, \ldots ,L_{m_1+m_2+2}$ such that $T$ acts trivially on these bundles after restriction to $Y_2$. We continue in this way for the fixed point components $Y_3,\ldots ,Y_{l-1}$.

For the remaining $m- \sum _{i=1}^{l-1} (m_i+1)\geq m_l$ line bundles we fix a lift of the $T$-action such that their local weights at $Y_l$ are all zero.

Let $v_i\in H_T^2(M;\Z)$ be the equivariant extension of $x_i$ given by the equivariant first Chern class of $L_i$. Let $y_1:=v_1\cdot \ldots \cdot v_{m_1+1}, y_2:=v_{m_1+2}\cdot \ldots \cdot v_{m_1+m_2+2},\ldots ,\linebreak y_l:=v_k\cdot v_{k+1}\cdot \ldots \cdot v_m$, where $k:={\sum _{i=1}^{l-1} (m_i+1)+1}$.

Let $v:=v_1\cdot \ldots \cdot v_m=y_1\cdot \ldots \cdot y_l\in H^{2m}_T(M;\Z)$ be the product of the equivariant first Chern classes of the line bundles. By construction $v$ is an equivariant extension of $x_1\cdot \ldots \cdot x_m$.

Also by construction the restriction $j_{Y_i}^*(y_i)$ of $y_i$ to $Y_i$ vanishes for $i=1,\ldots ,l-1$. If $m+1>\sum _{i=1}^l (m_i+1)$ then $j_{Y_l}^*(y_l)$  also vanishes.  If $m+1=\sum _{i=1}^l (m_i+1)$ then $j_{Y_l}^*(y_l)=j_{Y_l}^*(x_k\cdot \ldots \cdot x_m)\in H^{2m_l}(Y_l;\Z)$.

By Fixed point formula \ref{fpf} we have
$$\int _M (x_1\cdot \ldots \cdot x_m)=\int _M  v=\sum _{i=1}^l a(Y_i),$$
where the local datum $a(Y_i)$ is equal to $\int _{Y_i} \frac {j_{Y_i}^*(v)} {e_T(\nu _{Y_i})}$. Note that $j_{Y_i}^*(y_i)=0$ implies $a(Y_i)=0$. Since $\int _M (x_1\cdot \ldots \cdot x_m)\neq 0$ we conclude that $a(Y_i)=0$ for $i<l$, $a(Y_l)\neq 0$ and $m+1=\sum _{i=1}^l (m_i+1)$. Since the local datum
$$a(Y_l)=\int _{Y_l}j_{Y_l}^*(y_l)\cdot \frac {j_{Y_l}^*(y_1\cdot \ldots \cdot y_{l-1})} {e_T(\nu _{Y_i})}$$ is non-zero $j_{Y_l}^*(y_l)$ does not vanish.

Now $j_{Y_l}^*(y_l)\in H^{2m_l}(Y_l; \Z )$ is equal to $\tilde x_k\cdot \ldots \cdot \tilde x_m\in H^{2m_l}(Y_l; \Z )$, where $\tilde x_i\in H^2(Y_l;\Z)$ is the restriction of $x_i$ to $Y_l$. Hence, $Y_l$ is $H^2$-orientable. Of course an analogous reasoning applies to any $Y_i$. This finishes the proof.\proofend

\begin{corollary} Let $M$ be an oriented closed connected manifold of dimensional $2m$ with $b_1(M)=0$. Let $T$ be a torus which acts smoothly on $M$ with fixed point components $Y_1,\ldots ,Y_l$ of dimension $\dim Y_i=2m_i$.

Suppose $M$ is $H^2$-orientable. Then $m+1\leq \sum _{i=1}^l (m_i+1)$. Moreover, if equality holds then each component $Y_i$ is $H^2$-orientable.\proofend
\end{corollary}

In the corollary the condition $b_1(M)=0$ is used to guarantee that the action lifts to the line bundles. If $b_1(M)\neq 0$ then this may still be true. For example, if the torus action extends to an action by a semi-simple Lie group $G$ then the $G$-action and hence the torus action lifts to any line bundle \cite{St61,Su63,HaYo76}.

Note that by the theorem above the number of (not necessarily isolated) fixed points of a torus action on a $2m$-dimensional manifold which is equivariantly $H^2$-orientable is $\geq m+1$. If the condition $m+1= \sum _{i=1}^l(m_i+1)$ holds then the theorem implies that the sum of even Betti numbers of $M^T$ is $\geq m+1$. Both lower bounds can be realized by linear actions on $\C P^m$.

We close this section with some remarks concerning actions with few fixed points. Let $M$ be an oriented closed manifold of positive dimension. Assume a torus $T$ acts on $M$ smoothly with at least one fixed point.

The first thing we would like to recall is that the action cannot have precisely one fixed point. To see this let us apply Fixed point formula \ref{fpf} to the unit element in the integral cohomology ring of $M$ which extends trivially to the equivariant class $v=1\in H^0_T(M;\Z)$.

Since $M$ is of positive dimension $\int _M v=0$. The local datum at a fixed point component $Y$ is equal to $\int _{Y} \frac {j_{Y}^*(v)} {e_T(\nu _{Y})}$. It is easy to verify that the local datum is non-zero if $Y$ is a point. Hence, $T$ cannot act with precisely one fixed point. A proof using index theory which also applies to the action by a group of order $p^l$, $p$ an odd prime, was given by Atiyah and Bott in \cite[Thm. 7.1]{AtBoII68}, for a proof using bordism see \cite[\S 8]{CoFl66}, for a cohomological proof see \cite[IV Cor. 2.3]{Br72}.

If the torus $T$ acts on $M$ with precisely two fixed points then the representations at the two fixed points are equivalent \cite[6. Appendix]{Ko84}. This can be shown using the rigidity of the signature following the line of reasoning in \cite[Thm. 7.15]{AtBoII68}.

For actions with precisely three fixed points the local situation at the fixed points is more complicated since there are many more examples of such actions including homogeneous actions on the complex projective space $\C P^2$, the quaternionic projective plane $\H P^2$, the rational quaternionic projective plane $G_2/SO(4)$ and the Cayley plane $F_4/Spin(9)$. It would be interesting to understand to what extend the fixed point formulas and the rigidity and vanishing theorems of the following sections determine the local geometry and the bordism type for actions with three fixed points.

\section{Rigidity and vanishing of classical operators}\label{classicaloperatorsection}
In this section we review the rigidity and vanishing theorems for the signature operator, the Dirac operator and the Dolbeault operator. In particular cases the indices of these operators are obstructions to smooth circle actions and/or to actions with a lower bound on the codimension of the fixed point manifold. For example, by the celebrated $\hat A$-vanishing theorem of Atiyah and Hirzebruch \cite{AtHi70} the index of the Dirac operator obstructs smooth non-trivial circle actions on spin manifolds. And the signature obstructs $S^1$-actions for which the codimension of the fixed point manifold \wrt \ the involution in $S^1$ is larger than half of the dimension. More refined results based on the rigidity and vanishing theorems for elliptic genera will be described in the following section. Towards the end of this section we will also mention some relations to manifolds of positive or nonnegative curvature.

Throughout the section we will use the cohomological version of the index theorem and the Lefschetz fixed point formula which will be reviewed first.

\subsection{Cohomological form of the Lefschetz fixed point formula}\label{index theorem in cohomology}
As explained by Atiyah and Singer the $K$-theoretical Lefschetz fixed point formula \ref{fpfindK} can be translated into cohomology using the Chern character. We will only illustrate this in the basic examples below. For the details we refer to \cite[\S 3, Lefschetz Theorem 3.9]{AtSiIII68} (see also \cite[Ch. 5]{HiBeJu92}).

Let $M$ be an oriented closed Riemannian $2m$-dimensional manifold, $E_0,\ldots ,E_l$ complex vector bundles over $M$ and $D_i:\Gamma (E_i)\to \Gamma (E_{i+1})$ differential operators such that
$$\ldots \overset {D_{i-1}}\longrightarrow \Gamma (E_i) \overset {D_{i}}\longrightarrow \Gamma (E_{i+1})\overset {D_{i+1}}\longrightarrow \Gamma (E_{i+2})\to \ldots $$
is an elliptic complex $D$. Suppose the bundles $E_i$ are associated to the tangent bundle. Then the index $\ind (D)$ of the complex $D$ can be computed cohomologically by the formula
$$\ind (D)=\int _M  \frac {\sum _{i=0}^l (-1)^i\cdot ch(E_i)}{e(TM)}\cdot td(TM_\C ).$$
Here $td(TM_\C )=\prod _{i=1}^m\left (\frac {x_j}{1-e^{-x_j}}\cdot \frac {-x_j}{1-e^{x_j}}\right)$ is the Todd class and $\pm x_1,\ldots ,\pm x_m$ are the roots of the complexified tangent bundle $TM_{\C }$. The same is true if the bundles $E_i$ are associated to a reduction of the structure group of the stable tangent bundle satisfying a certain condition (see Condition 2.17 in \cite{AtSiIII68}) which will be always satisfied in the situations considered in this paper.

Next suppose a compact Lie group $G$ acts smoothly on $M$ and the complex $D$ is a $G$-equivariant elliptic complex. Let $g\in G$. Recall from Lefschetz fixed point formula \ref{fpfindK} that the index $\ind (D)\in R(G)$ of $D$ evaluated at $g$ is given by a sum of local data
$$\ind (D)(g)=\sum _{Y\subset M^g} a_Y.$$
In terms of cohomology the local datum $a_Y$ at a connected fixed point component $Y$ can be computed by the following recipe (see \cite[\S 3]{AtSiIII68}, \cite[Ch. 5.6]{HiBeJu92} for details).

Let us first replace $G$ by the compact Lie group generated by $g$. The normal bundle of $Y$ splits \wrt \ the action of $G$ as a sum of subbundles corres\-ponding to the different representations. The roots of each subbundle refine to equivariant roots. Similarly the roots of the bundles $E_i$ refine to equivariant roots \wrt \ the action of  $G$. Now consider the cohomological expression
\begin{equation}\label{cohomology expression}\frac {\sum _{i=0}^l (-1)^i\cdot ch(E_i)}{e(TM)}\cdot \prod _{i=1}^m\left (\frac {x_j}{1-e^{-x_j}}\cdot \frac {-x_j}{1-e^{x_j}}\right)\end{equation}
in the index formula. Note that (\ref{cohomology expression}) is an expression in the roots of $M$ and the roots of the bundles $E_i$. At a fixed point component $Y$, as indicated above, the roots of $M$ refine to the equivariant roots of the normal bundle and the roots of $Y$.

To obtain the cohomological expression for the local datum $a_Y$ one replaces in expression (\ref{cohomology expression}) the roots of $M$ by its equivariant roots evaluated at $g$, replaces the roots of the bundles $E_i$ by its equivariant roots evaluated at $g$, replaces the Euler class $e(TM)$ by $e(TY)$ and integrates over $Y$. We will make this explicit in the basic examples below. For details we refer to  \cite[\S 3]{AtSiIII68}, \cite[Ch. 5.6]{HiBeJu92}.

\subsection{The signature operator}
Let $M$ be an oriented closed Riemannian $2m$-dimensional manifold. Then one can define the signature operator on $M$ acting on certain spaces of differential forms \cite[\S 6]{AtSiIII68}. Its index can be identified via Hodge theory with the signature $sign(M)\in \Z $. The latter is a purely cohomological invariant which is defined as the signature of the intersection form of $M$. By the index theorem the index of the signature operator is equal to the $L$-genus. This gives an index theoretical proof of the celebrated Hirzebruch signature theorem.

Suppose a compact Lie group $G$ acts on $M$ by isometries. Then the signature operator refines to a $G$-equivariant operator, the equivariant signature operator. Its index is a virtual complex $G$-representation and is denoted by $sign(M)_G\in R(G)$. The latter can be described in terms of the $G$-action on certain subspaces of the space of harmonic forms \cite[p. 578]{AtSiIII68}. If $G$ is connected then the action of $G$ on harmonic forms is trivial, by homotopy invariance, and one gets the following

\begin{theoremsignrigid} If $G$ is connected then the equivariant signature $sign(M)_G$ is constant as a character of $G$, i.e. $sign(M)_G=sign(M)\in \Z$.\proofend
\end{theoremsignrigid}
Because of the significance of rigidity phenomenons for the subsequent discussion we will indicate a different argument using the cohomological form of Lefschetz fixed point formula \ref{fpfindK}.

The cohomological expression (\ref{cohomology expression}) for the index of the signature operator simplifies (see \cite[(6.4)]{AtSiIII68}, \cite[p. 65]{HiBeJu92}) and one obtains
$$sign (M) =\int _M e(TM)\cdot \prod_{i=1}^{2m}\frac  {1+e^{-x_i}}{1-e^{-x_i}}.$$
Suppose that $G=S^1$ is a circle acting isometrically on $M$ and $Y\subset M^{S^1}$ is a fixed point component. We will consider the complex structure on the normal bundle $\nu _Y$ of $Y\subset M$ which is induced by the action of a generator $\lambda =e^{2\pi i\cdot z}\in S^1$, $z>0$ small, on the fibers of $\nu _Y$. The complex vector bundle $\nu _Y$ splits as a sum of complex subbundles corresponding to the irreducible representations of $S^1$ of positive $S^1$-weight. We equip $Y$ with the orientation compatible with the complex structure of $\nu _Y$ and the orientation of $M$.

The equivariant roots of the $S^1$-equivariant complex normal bundle $\nu _Y$ are of the form $x_i+2\pi i \cdot m_i\cdot z$, where all $m_i$ are positive integers. Let $\pm y_j$ denote the roots of $Y$, i.e. the equivariant roots with zero $S^1$-weight.

Using the recipe above one computes for the local datum $a_Y(\lambda )$, $\lambda \in S^1$, that
\begin{equation}\label{signlocaldatum}a(Y)(\lambda)=\int _Y e(TY)\cdot \prod_{j}\frac  {1+e^{-y_j}}{1-e^{-y_j}}\cdot \prod_{i}\frac  {1+e^{-(x_i+2\pi i \cdot m_{i}\cdot z)}}{1-e^{-(x_i+2\pi i \cdot m_{i}\cdot z)}}$$
$$=\int _Y e(TY)\cdot \prod_{j}\frac  {1+e^{-y_j}}{1-e^{-y_j}}\cdot \prod_{i}\frac  {1+e^{-x_i}\cdot \lambda ^{-m_{i}}}{1-e^{-x_i}\cdot \lambda ^{-m_{i}}}.\end{equation}

Note that each local datum $a_Y(\lambda )$, $\lambda \in S^1$, extends to a meromorphic function on the complex plane with possible poles only on $S^1$ and which is bounded in $\infty $. On the other hand the index is a character of $S^1$ and, hence, determines a meromorphic function on the complex plane with possible poles only in $0$ and $\infty$. Now both meromorphic functions, the character and the sum of local data, agree on a dense set of $S^1$. Hence, the two functions are equal and bounded holomorphic on $\C $. It follows that they are constant in $\lambda $ and $sign(M)_{S^1}=sign(M)\in \Z$ (for details see for example \cite{BoTa89,HiBeJu92}).

For a connected compact Lie group $G$ acting isometrically on $M$ the rigidity of $sign(M)_G$ follows since a character of $G$ is constant if and only if its restriction to all $S^1$-subgroups is constant.

\begin{remark} Taking the limit $\lambda \to \infty$ in (\ref{signlocaldatum}) it follows from the rigidity of the signature that
$$sign (M)=\sum _{Y\subset M^{S^1}} sign (Y) .$$
\end{remark}

The rigidity is not true in general if $G$ is non-connected. Nevertheless there are interesting applications of the fixed point formula for such actions, see \cite{AtBoII68,AtSiIII68,HiZa74}. Related to the rigidity of the signature is the following vanishing result.

\begin{theoremsignvanish}\label{theoremsignvanish} Suppose $S^1$ acts smoothly on an oriented closed manifold $M$. Let $\sigma \in S^1$ be of order two and let $M^\sigma$ be its fixed point manifold. If $\dim M^\sigma < \frac 1 2 \dim M$ then $sign(M)=0$.\proofend
\end{theoremsignvanish}
We note that after averaging a Riemannian metric over $S^1$ we may assume that $S^1$ acts by isometries. The vanishing result \ref{theoremsignvanish} is due to Hirzebruch \cite[p.~153]{Hi68}. He used the cohomological version of the Lefschetz fixed point formula for the equivariant signature \wrt \ the action of the involution $\sigma \in S^1$ to show that
\begin{equation}\label{signselfintersection} sign(M)_{S^1}(\sigma )=\sum _{Y\subset M^\sigma }sign (Y\circ Y),\end{equation}
where $Y\circ Y$ denotes a transversal self-intersection in $M$ of the fixed point component $Y$ (for an elementary proof not relying on the index theorem see \cite{JaOs69}). If $\dim Y<\frac 1 2 \dim M$ then the transversal self-intersection is empty. By the rigidity of the signature $sign(M)=sign(M)_{S^1}(\sigma )$. Hence, one obtains
$$\dim M^\sigma < \frac 1 2 \dim M \implies sign (M)=0.$$
This {\em rigidity implies vanishing} phenomenon will be a recurring theme in the following sections.

By the vanishing result \ref{theoremsignvanish} the non-vanishing of the signature implies the existence of a fixed point component $Y\subset M^\sigma$ with $2\dim Y \geq \dim M$ (see \cite[\S 27]{CoFl64} for an earlier result involving the Euler characteristic). In this connection we would also like to mention that by a result of Boardman \cite{Boa67} a manifold $M$ bounds as an unoriented manifold if it admits a smooth action by an involution $\sigma $ with $\dim M^\sigma < \frac 2 5\dim M$ (see also \cite{KoSt78}).

\subsection{The Dirac operator}

Let $M$ be a closed $4m$-dimensional spin manifold. Then one can define a Dirac operator acting on spinors of $M$ \cite[\S 5]{AtSiIII68}. By the index theorem the index of the Dirac operator is the $\hat A$-genus $\hat A(M)\in \Z$, a topological invariant which can be computed from the Pontrjagin classes of $M$.

Now suppose a compact Lie group $G$ acts on $M$ preserving the spin structure. Then the Dirac operator refines to an equivariant operator with equivariant index $\hat A(M)_G\in R(G)$. In \cite{AtHi70} Atiyah and Hirzebruch proved the following

\begin{theoremAdachvanish} Let $M$ be a connected spin manifold and let $G$ be a connected compact Lie group which acts non-trivially on $M$ preserving the spin structure. Then $\hat A(M)_G=0\in R(G)$.\proofend
\end{theoremAdachvanish}

The proof of Atiyah and Hirzebruch is based on the cohomological form of Lefschetz fixed point formula \ref{fpfindK} applied to the symbol of the equivariant Dirac operator. The argument roughly runs as follows.

Suppose first that $G=S^1$ is a circle and $Y\subset M^{S^1}$ denotes a fixed point component. By the cohomological version of the index theorem the index of the Dirac operator is equal to
$$\hat A(M)=\int _M e(TM)\cdot \prod _{i=1}^{2m}\frac  {1}{e^{x_i/2}-e^{-x_i/2}}.$$
With the notation of the last subsection the recipe gives for the local datum at $Y$:
\begin{equation}\label{Ahatlocaldatum}a(Y)(\lambda)=\int _Y e(TY)\cdot \prod_{j}\frac  {1}{e^{y_j/2}-e^{-y_j/2}}\cdot \prod_{i}\frac  {1}{e^{(x_i+2\pi i \cdot m_{i}\cdot z)/2}-e^{-(x_i+2\pi i \cdot m_{i}\cdot z)/2}}$$
$$=\int _Y e(TY)\cdot \prod_{j}\frac  {1}{e^{y_j/2}-e^{-y_j/2}}\cdot \prod_{i}\frac  {1}{e^{x_i/2}\cdot \lambda ^{m_{i}/2}-e^{-x_i/2}\cdot \lambda ^{-m_{i}/2}}.\end{equation}
One computes that each local datum $a_Y(\lambda )$, $\lambda \in S^1$, extends to a meromorphic function on the complex plane which vanishes in $0$ and $\infty$ and with possible poles only on $S^1$. On the other hand the equivariant index of the Dirac operator is a character of $S^1$ and, hence, determines a meromorphic function on the complex plane with possible poles only in $0$ and $\infty$. Now both meromorphic functions, the character and the sum of local data, agree on a dense set of $S^1$. Hence, the two functions are equal, bounded holomorphic on $\C $ and vanish in $\infty$. It follows that they are constant zero in $\lambda $ and $\hat A(M)_{S^1}=\hat A(M)=0\in \Z$ (for details see \cite{AtHi70}, see also \cite{BoTa89,HiBeJu92}).

For a connected compact Lie group $G$ which acts non-trivially on $M$ preserving the spin structure the equivariant index $\hat A(M)_G$ also vanishes. This follows from the above since $\hat A(M)_{S^1}$ is rigid for any subgroup $S^1\subset G$, $\hat A(M)_{S^1}$ vanishes if $S^1$ acts non-trivially and a character of $G$ vanishes if and only if its restriction to all $S^1$-subgroups vanishes.

\subsection{The Dolbeault complex}

In this subsection we assume that $M$ is an almost complex closed manifold of real dimension $2m$. The results below extend to stably almost complex manifolds. To simplify the discussion we will restrict to the unstable case.

Using the Dolbeault operator $\bar \partial $ one can define the Dolbeault complex
$$ \ldots \overset {\bar \partial }\longrightarrow \Gamma (A^{p,q-1})\overset {\bar \partial }\longrightarrow \Gamma (A^{p,q})\overset {\bar \partial }\longrightarrow \Gamma (A^{p,q+1})\overset {\bar \partial }\longrightarrow \ldots ,$$
where $\Gamma (A^{p,q})$ is the space of sections in the vector bundle $\Lambda ^p(T^*M)\otimes \Lambda ^q(\overline {T^*M})$.

For fixed $p$ let $\chi ^p(M)$ denote the index of this complex (we are following the notation of \cite{HiBeJu92}). By the index theorem (see \cite[\S 4]{AtSiIII68}) $\chi ^p(M)$ can be computed from the Chern classes of $M$. In particular, $\chi ^0(M)$ is equal to the Todd genus $Td (M)$. For $\chi_y(M):=\sum _{p=0}^m \chi ^p(M)\cdot y^p\in \Z [y]$ one obtains by the index theorem that
$$\chi _y(M)=Td_y(M):=\int _M \prod _{i=1}^m \left ( (1+y\cdot e^{-x_i})\cdot \frac {x_i}{1-e^{-x_i}}\right ). $$
Here $x_1,\ldots ,x_m$ are the roots of the almost complex structure.
This gives a generalization of the Hirzebruch-Riemann-Roch theorem to almost complex manifolds.

Now suppose a compact Lie group $G$ acts on $M$ preserving the almost complex structure. Then the Dolbeault complex refines to an equivariant complex and $\chi ^0(M)=Td(M)\in \Z $, resp. $\chi _y(M)=Td_y(M)\in \Z[y]$, refine to equivariant indices which we denote by $Td (M)_G\in R(G)$, resp. $Td_y (M)_G\in R(G)[y]$. If $G$ is connected then, like for the signature and Dirac operator, one has rigidity.

\begin{theoremToddrigid}\label{theoremToddrigid} Let $M$ be an almost complex manifold and $G$ a compact connected Lie group which acts on $M$ preserving the almost complex structure. Then the equivariant index $Td _y(M)_G\in R(G)[y]$ is rigid, i.e. $Td_y(M)_G=Td _y(M)\in \Z[y]$. In particular, $Td (M)_G$ is rigid.\proofend
\end{theoremToddrigid}
This was first shown by Lusztig \cite{Lu71} for holomorphic actions with isolated fixed points using the Lefschetz fixed point formula (see also Kosniowski \cite[Thm. 5]{Ko70}). The proof carries directly over to almost complex manifolds, for details see for example \cite[\S 5.7]{HiBeJu92}. The genus $Td _y$ extends to a two-parameter genus $Td _{x,y}$ which is rigid for compact connected Lie group actions on (stably) almost complex manifolds preserving the structure and is characterized by this property, as shown by Krichever \cite{Kr74} and Musin \cite{Mu11}, respectively.

The rigidity of $Td_y$ generalizes to the rigidity of elliptic genera of level $N$. Whereas the rigidity of $Td_y$ does not lead directly to a vanishing theorem the rigidity of elliptic genera of level $N$ implies vanishing results which are similar to the $\hat A$-vanishing theorem above (see Section \ref{levelNsection}).

\subsection{Connections to curvature}

We close this section by mentioning some relations between the indices of classical operators and the existence of metrics of nonnegative or positive curvature.

Concerning the signature one has the following absolute bound which is a consequence of Gromov's Betti number theorem \cite{Gr81}.

\begin{theorem}[Gromov] If $M$ is connected and of nonnegative
sectional curvature then the sum of the Betti numbers of $M$ is
bounded from above by a constant $C$ which depends only on the
dimension of $M$. In particular, $\vert
sign(M)\vert <C$.\proofend \end{theorem}

The Dirac operator is intimately related to scalar curvature by the work of Licherowicz, Gromow-Lawson, Schoen-Yau, Stolz and others. In \cite{Li63} Lichnero\-wicz used the Weitzenb\"ock-formula for the Dirac operator to prove

\begin{theorem}[Lichnerowicz] If $M$ is a spin manifold
with positive scalar curvature then $\hat A(M)=0$.\proofend \end{theorem}

For K\"ahler manifolds of positive Ricci curvature one has Bochner's theorem \cite{Bo46}, \cite[\S IV, Cor. 11.12]{LaMi89}.

\begin{theorem}[Bochner] If $M$ is a connected K\"ahler manifold of positive Ricci curvature then $Td(M)=1$.\proofend \end{theorem}

\section{Elliptic genera of level $2$}\label{level2section}
In this section and the following section we review the rigidity and vanishing theorems for elliptic genera. The vanishing results are consequences of the rigidity and are related to finite cyclic subactions for which the fixed point manifold is of {\em large} codimension. This gives new information on the fixed point manifold. We will also remark on connections to manifolds of positive sectional curvature.

We begin the discussion with the elliptic genus of level $2$ of a spin manifold $M$ which, according to Witten, should be thought of as the equivariant signature of the free loop space ${\mathcal L}M$. This genus emerged from the work of Landweber, Ochanine, Stong, Witten and others in the 1980s (see \cite{La88} and references therein).

\subsection{Twisted indices} To define the elliptic genus of level $2$ we need to recall the computation of twisted signature and twisted Dirac operators via the index theorem.

Let $M$ be an oriented closed $4m$-dimensional Riemannian manifold and $E\to M$ a complex vector bundle. After choosing a connection on $E$ one can define the signature operator twisted with $E$. Its index will be denoted by $sign(M,E)\in \Z $ and is called a {\em twisted signature}. By the index theorem \cite{AtSiIII68} $sign(M,E)$ can be computed from the Pontrjagin classes of $M$ and the Chern character of $E$:
$$sign(M,E) = \int _M e(TM)\cdot \prod_{i=1}^{2m}\frac  {1+e^{-x_i}}{1-e^{-x_i}}\cdot ch(E)$$

Now suppose a compact connected Lie group $G$ acts on $M$ by isometries and the action lifts to the bundle $E$. Then the twisted signature operator refines to an equivariant operator and the twisted signature refines to a virtual complex $G$-representation denoted by $sign(M,E)_G\in R(G)$.

In contrast to the ordinary signature a twisted signature does not need to be rigid. For example if $M$ is the projective plane $\C P^2$ and $E$ is its complexified tangent bundle $TM_\C $ then for any nontrivial homogeneous action of $S^1$ on $M$ the equivariant twisted signature $sign(M,E)_{S^1}$ is a non-trivial character.

Next let $M$ be spin and let $E\to  M$ be a complex vector bundle. After choosing a connection on $E$ one can define the Dirac operator twisted with $E$. Its index will be denoted by $\hat A(M,E)\in \Z $ and is called a {\em twisted Dirac index}. By the index theorem \cite{AtSiIII68} $\hat A(M,E)$ can be computed from the Pontrjagin classes of $M$ and the Chern character of $E$:
$$\hat A(M,E)=\int _M e(TM)\cdot \prod _{i=1}^{2m}\frac  {1}{e^{x_i/2}-e^{-x_i/2}}\cdot ch(E)$$

Now suppose a compact connected Lie group $G$ acts on $M$ by isometries and the action lifts to the spin structure and to the bundle $E$. Then the twisted Dirac operator refines to an equivariant operator and the twisted Dirac index refines to a virtual complex $G$-representation denoted by $\hat A(M,E)_G\in R(G)$. Again there are examples where $\hat A(M,E)_G$ is non-zero and not rigid.

\subsection{Signature of the free loop space (elliptic genus of level $\boldsymbol 2$)} Following Witten the elliptic genus of level $2$ should be thought of (up to a normalization factor) as the equivariant signature of the free loop space ${\mathcal L}M$. Here the action on ${\mathcal L}M$ is the natural action of $S^1$ by rotation of loops. The fixed point set of this action is the manifold of constant loops $M=({\mathcal {L}}M)^{S^1}$.

Until now a solid mathematical description for Witten's heuristic is still missing. However, by applying the Lefschetz fixed point formula formally to the natural $S^1$-action on ${\mathcal L}M$ Witten derived the following well-defined invariant for the underlying manifold $M$ (see \cite{Wi86}, see also \cite{HiBeJu92}).

\begin{definition}\label{loop space signature}
$$sign({\mathcal {L}}M):=sign(M,\bigotimes_{n=1}^\infty S_{q^n} TM_\C\otimes \bigotimes_{n=1}^\infty \Lambda _{q^n}TM_\C)$$

$$=\sum _{m=0}^\infty sign(M,R_m)\cdot q^m=sign(M)+2sign(M, TM_\C )\cdot q + \ldots .$$
\end{definition}
Here
$S_t:=\sum _iS^i\cdot t^i$ (resp. $\Lambda _t:=\sum _i \Lambda
^i\cdot t^i$) denotes the symmetric (resp. exterior) power
operation and $TM_\C$ is the complexified tangent bundle. The bundles $R_m$ are computed by taking the coefficients in the $q$-power series
$$\bigotimes_{n=1}^\infty S_{q^n} TM_\C\otimes \bigotimes_{n=1}^\infty \Lambda _{q^n}TM_\C,$$ the first few are
$$R_0=1, R_1=2\, TM_\C , R_2=2\, ( TM_\C +TM_\C\otimes TM_\C ),\ldots .$$
We shall call $sign({\mathcal {L}}M)$ the {\em signature of the free loop space} of $M$. The definition does not require a spin structure. By the index theorem $sign({\mathcal {L}}M)$ can be computed from the Pontrjagin classes of $M$.

The series $sign({\mathcal {L}}M)$ converges for $\vert q\vert <1$ to a meromorphic function $\Phi (M)(\tau)$ on the upper half plane, where $q=e^{2\pi i \tau}$. It turns out that $\Phi (M)(\tau)$ is a modular function of weight zero for the congruence subgroup
$$\Gamma _1(2):=\{A=(\begin{smallmatrix}a & b\\ c & d
\end{smallmatrix}) \in
SL_2(\Z)\mid (\begin{smallmatrix}a & b\\ c & d
\end{smallmatrix})\equiv (\begin{smallmatrix}1 & b\\ 0 & 1
\end{smallmatrix}) \bmod 2\},$$
 (for an explanation of these facts see \cite{Wi86,BoTa89,HiBeJu92}).

It follows from the construction that $\Phi $ is a ring homomorphism from the oriented bordism ring to the ring of modular functions for $\Gamma _1(2)$. Up to a normalization constant $\Phi $ is equal to the universal elliptic genus of level $2$ considered by Ochanine, Landweber and Stong (see \cite{La88}). We will also call $\Phi $ an {\em elliptic genus (of level $2$)}.

By the definition of $\Phi (M)$ the expansion of this modular function in the cusp $i\infty$ ({\em signature cusp}) is equal to $sign({\mathcal {L}}M)$. Another interesting expansion is obtained by taking the cusp $0$. The expansion of $\Phi (M)$ in the cusp $0$ (the {\em Dirac cusp}) can be computed to be a Laurent series in $q^{\frac 1 2}$ with a pole of finite order.

Now assume $M$ is spin. After replacing $q^{\frac 1 2}$ by $q$ this Laurent series can be shown, with the help of the index theorem, to be equal to the following series of twisted Dirac indices

$$C\cdot \hat
A(M,\bigotimes _{n=2m+1>0}\Lambda _{-q^n}TM_\C \otimes \bigotimes
_{n=2m>0}S_{q^n}TM_\C)$$

$$=C\cdot (\hat A(M)
-\hat A(M,TM_\C)\cdot q +\hat A(M,\Lambda ^2{TM_\C}+TM_\C)\cdot q^2+\ldots ).$$

Here $C\in \Z[q^{-1}][[q]]$ is a constant depending only on the dimension of $M$ which has a pole of order $\frac {\dim M}8$.

Thus the elliptic genus of level $2$ provides a connection between twisted signatures and twisted Dirac indices. This has led to many applications including an elegant proof of the divisibility theorem for the signature of Ochanine-Rochlin \cite[Ch. 8]{HiBeJu92}.

\subsection{Rigidity}

Let $M$ be a closed spin manifold and let $G$ be a compact Lie group which acts on $M$ preserving the spin structure. Then $G$ acts on the spin principal bundle and on all vector bundles associated to the spin principal bundle. In particular, the coefficients in the expansion of $\Phi (M)$ in the signature cusp (resp. Dirac cusp) refine to equivariant twisted signatures (resp. equivariant Dirac indices). That is, the expansion of the equivariant elliptic genus $\Phi(M)_G$ in the signature cusp takes the form

$$sign({\mathcal {L}}M)_G=sign(M,\bigotimes_{n=1}^\infty S_{q^n}
TM_\C\otimes \bigotimes_{n=1}^\infty \Lambda _{q^n}TM_\C)_G$$
$$=sign(M)_G+2sign(M, TM_\C )_G\cdot q + \ldots\in R(G) [[q]]$$
and the expansion in the Dirac cusp takes the form
$$C\cdot \hat
A(M,\bigotimes _{n=2m+1>0}\Lambda _{-q^n}TM_\C \otimes \bigotimes
_{n=2m>0}S_{q^n}TM_\C)_G$$
$$=C\cdot (\hat A(M)_G
-\hat A(M,TM_\C)_G\cdot q +\hat A(M,\Lambda ^2{TM_\C}+TM_\C)_G\cdot q^2+\ldots )\in R(G) [[q]].$$

We are now in the position to state the
\begin{rigiditytheorem}[Witten, Taubes, Bott-Taubes] Let $M$ be a closed spin manifold and let $G$ be a compact connected Lie group which acts on $M$ preserving the spin structure. Then $\Phi(M)_G$ and its expansion in the two cusps are rigid, i.e. each coefficient in the expansions above is constant as a character of $G$.\proofend\end{rigiditytheorem}

The rigidity of the elliptic genus was conjectured by Witten \cite{Wi86} and then shown by Taubes \cite{Ta89} and Bott-Taubes \cite{BoTa89}. The proof of Bott and Taubes uses the Lefschetz fixed point formula together with arguments from elliptic function theory and a technical transfer argument. A simpler proof based on the Lefschetz fixed point formula and modularity invariance of Jacobi functions was later given by Liu \cite{Li95}.

\subsection{Higher order vanishing} The rigidity of the elliptic genus leads to vanishing results for finite cyclic subactions. For involutions such results were shown by Hirzebruch and Slodowy \cite{HiSl90} using a formula involving the transversal self-intersection of the fixed point manifold similar to formula (\ref{signselfintersection}) for the signature. An extension to finite cyclic subactions of higher order was given in \cite{De01b,De07} using different arguments.

Let $M$ be a closed spin manifold with an action of $S^1$ by isometries. If the $S^1$-action does not lift to the spin structure then the non-trivial two-fold covering of $S^1$ lifts to the spin structure. In this case one can show, using the Lefschetz fixed point formula and the rigidity theorem, that the equivariant elliptic genus \wrt \ the two-fold $S^1$-action vanishes identically. This implies the vanishing of $sign({\mathcal {L}}M)_{S^1}$ for the original $S^1$-action. Hence, for the following vanishing results we may, without mentioning, always assume that the action lifts to the spin structure.

Let $\sigma \in S^1$ be of order two. Consider the fixed point manifold $M^\sigma$ and its transversal self-intersection $M^\sigma \circ M^\sigma $. In \cite{HiSl90} Hirzebruch and Slodowy used the cohomological version of the Lefschetz fixed point formula \ref{fpfindK} to compute that $sign({\mathcal {L}}M)_{S^1}(\sigma)$ is equal to $sign({\mathcal {L}}(M\circ M))$. Note that this is completely analogous to formula (\ref{signselfintersection}) for the ordinary signature. Now by the rigidity it follows that
$$sign({\mathcal {L}}M)=sign({\mathcal {L}}(M\circ M)) \text{ and }\Phi (M)=\Phi(M\circ M).$$ Changing to the Dirac cusp and comparing the pole order Hirzebruch and Slodowy \cite{HiSl90} arrive at the following vanishing result.

\begin{theorem}[Hirzebruch-Slodowy]\label{HiSlvanishing} Let $M$ be a closed spin manifold and let $S^1$ act smoothly on $M$. Let $\sigma \in S^1$ be of order two.
If the codimension of $M^\sigma$ is $ > 4r$
 then the first $(r+1)$ coefficients in the expansion $$\hat
A(M,\bigotimes _{n=2m+1>0}\Lambda _{-q^n}TM_\C \otimes \bigotimes
_{n=2m>0}S_{q^n}TM_\C)$$ vanish.\proofend\end{theorem}
Here the codimension of a submanifold is defined as the minimum of the codimension of its connected components.

Note that this result provides a new way to prove the Atiyah-Hirzebruch vanishing for the $\hat A$-genus. In fact, if $S^1$ acts effectively and $M$ is connected then $\codim M^\sigma>0$ and, hence, $\hat A(M)=0$.

If $\codim M^\sigma>4$ then $\hat A(M)=\hat
A(M,TM_\C)=0$, if $\codim M^\sigma>8$ then $\hat A(M)=\hat
A(M,TM_\C)=\hat A(M,\Lambda ^2{TM_\C})=0$ and so on. Thus, if the codimension of $M^\sigma$ increases one obtains {\em  higher order vanishing} for the elliptic genus (see \cite{HiSl90}).

This higher order vanishing phenomenon generalizes to elements $\sigma \in S^1$ of arbitrary finite order. Suppose $S^1$ acts effectively on $M$ and let $\sigma \in S^1$ be of order $o\geq 2$. Then one can show \cite{De01a,De07} the following

\begin{theorem}\label{highervanishing} If $\codim M^\sigma
> 2 o \cdot r$ then the first $(r+1)$ coefficients in the expansion $$\hat
A(M,\bigotimes _{n=2m+1>0}\Lambda _{-q^n}TM_\C \otimes \bigotimes
_{n=2m>0}S_{q^n}TM_\C)$$ vanish.\proofend\end{theorem}

The proof of Hirzebruch-Slodowy only works for $\sigma$ of order two. The proof of Theorem \ref{highervanishing} involves the rigidity, arguments from function theory and a detailed study of the local contributions in the fixed point formula. It also leads to higher order vanishing for elliptic genera of level $N$ (see Section \ref{levelNsection}).

\begin{example}
If $\sigma $ is of order $3$ and $\codim M^\sigma>6$ then $\hat A(M)=\hat
A(M,TM_\C)=0$.\end{example}

\bigskip

\subsection{Connections to positive curvature}
In the last subsection we have described various situations where the existence of isometric actions for which the fixed point manifold has large codimension forces certain Dirac indices to vanish. On the other hand the fixed point components of isometric actions are totally geodesic submanifolds and these submanifolds enjoy in the presence of positive sectional curvature special properties.

For example Frankel \cite{Fr61} used a Synge type argument to prove that two totally geodesic connected submanifolds $N_1$, $N_2$ in a manifold $M$ of positive sectional curvature must intersect if $\dim N_1 +\dim N_2 \geq \dim M$. By combining Frankel's intersection theorem with the vanishing results for elliptic genera from the last subsection one obtains obstructions to the existence of metrics of positive sectional curvature with symmetry. More precisely such obstruction arise for closed positively curved spin manifolds of arbitrarily large dimension with a fixed lower bound on the symmetry rank (the {\em symmetry rank} of a Riemannian manifold is defined as the rank of its isometry group). Improved obstructions can be obtained if one also uses Wilking's connectivity lemma for totally geodesic submanifolds \cite[Thm. 1]{Wi03}.

This interplay between elliptic genera and positive curvature was first studied in \cite{De01a,De02}. For example one has the following result (see \cite[Th. 2.1]{De01a}, \cite[Th. 3]{De07}).

\begin{theorem} Let $M$ be a connected spin manifold of dimension $> 12r-4$. Suppose $M$ admits a metric of positive sectional curvature and of symmetry rank
$\geq 2r$. Then the first $(r+1)$ coefficients in the expansion $$\hat
A(M,\bigotimes _{n=2m+1>0}\Lambda _{-q^n}TM_\C \otimes \bigotimes
_{n=2m>0}S_{q^n}TM_\C)$$ vanish.\proofend \end{theorem}

The theorem gives new obstructions to positive sectional curvature if a torus of rank $\geq 2$ acts isometrically and effectively. In the case of an isometric effective action by $S^1$ one can show the following result (see \cite[Th. 1.2]{De01a}, \cite[Th. 1]{De07}).

\begin{theorem} Let $M$ be closed $2$-connected of dimension $> 8$. If $M$ admits a metric of positive sectional curvature and of symmetry rank $\geq 1$ then $\hat A(M)=\hat A(M,TM_\C)=0$.\proofend \end{theorem}

\begin{corollary} Let $M$ be a closed $2$-connected manifold of dimension $> 8$ with $\hat A(M,TM_\C)\neq 0$. Then any metric of positive sectional curvature must be quite unsymmetric.\proofend \end{corollary}

We like to stress that the symmetry assumptions in the results above are rather mild since for a given symmetry rank no upper bound on the dimension of the manifold $M$ is required! Similar results also hold for manifolds with positive $k$th Ricci curvature (see for example \cite[Prop. 19]{De07}).

Much more is known for manifolds of positive sectional curvature with {\em large} symmetry and we refer the interested reader to the survey articles \cite{Wi07,Zi07} and references therein. Here we will only mention the following results concerning the symmetry rank. In \cite{Wi03} Wilking obtained strong classification results for manifolds $M$ of positive sectional curvature and lower {\em linear} bound on the symmetry, e.g. symmetry rank $\geq \dim M /6$.

More recently Kennard \cite{Ke13} has found severe restrictions on the cohomology ring (in terms of periodicity)  for manifolds of positive sectional curvature for lower {\em logarithmic} bound and derived consequences for the Hopf conjecture. Kennard combines the connectivity lemma of Wilking with a topological result on manifolds with periodic cohomology to show, among other things, that the Euler characteristic is positive for a $4m$-dimensional manifold $M$ of positive sectional curvature if the symmetry rank in $> 2 \log m$. There is also interesting subsequent work by Amann and Kennard (see for example \cite{AmKe14}). For the elliptic genus Weisskopf has proved higher order vanishing using Kennard's work.

\begin{theorem}[\cite{We13}] Let $M$ be a connected spin manifold of large dimension. Suppose $M$ admits a metric of positive sectional curvature and of symmetry rank $s$. Then, M is rationally
$4$-periodic or the first $2^{s-3}$ coefficients in the expansion $$\hat
A(M,\bigotimes _{n=2m+1>0}\Lambda _{-q^n}TM_\C \otimes \bigotimes
_{n=2m>0}S_{q^n}TM_\C)$$ vanish.\proofend\end{theorem}

In the next section we will discuss rigidity and vanishing results for elliptic genera of higher level.

One can also define an elliptic genus of level one. The definition is due to Witten, who derived this genus from a heuristic argument involving a hypothetical Dirac operator on the free loop space. This genus, called the Witten genus, is expected to have interesting properties if the free loop space is spin. The latter condition is satisfied if the underlying manifold $M$ is string, i.e. if $M$ is spin and $\frac {p_1}2 (M)=0$. One interesting property is the rigidity and vanishing of the equivariant Witten genus under certain group actions (see \cite{Li95,De94,De00}). Another conjectural property, formulated by Stolz, is the vanishing of the Witten genus for string manifolds of positive Ricci curvature. For more information on this conjecture see \cite{St96} (see also \cite{De09}).

\section{Elliptic genera of level $N$}\label{levelNsection}
In this section we first review the rigidity theorems for elliptic genera of level $N\geq 2$. We will then state higher vanishing results for these genera which are analogous to the ones given in Theorem \ref{HiSlvanishing} and Theorem \ref{highervanishing}. As in the level $2$ case these higher vanishing results are consequences of the rigidity and are related to finite cyclic subactions for which the fixed point manifold is of large codimension. Details as well as applications will appear in \cite{De}.

To define the elliptic genus of level $N$ we need to recall the index theory for twisted Dolbeault operators.

\subsection{Twisted Dolbeault indices}
In this subsection we assume that $M$ is an almost complex closed manifold of real dimension $2m$. The results below extend to stably almost complex manifolds. To simplify the discussion we will restrict to the unstable case.

Let $E\to M$ be a complex vector bundle. For fixed $p$ one can twist the Dolbeault complex
$$ \ldots \overset {\bar \partial }\longrightarrow \Gamma (A^{p,q-1})\overset {\bar \partial }\longrightarrow \Gamma (A^{p,q})\overset {\bar \partial }\longrightarrow \Gamma (A^{p,q+1})\overset {\bar \partial }\longrightarrow \ldots $$
with $E$ to obtain a new elliptic complex. Its index will be denoted by $\chi ^p(M,E)$. By the index theorem (see \cite[\S 4]{AtSiIII68}) $\chi ^p(M,E)$ can be computed from the Chern classes of $M$ and $E$. In particular, $\chi ^0(M,E)$ is equal to
$$Td (M,E):=\int _M \left (\prod _{i=1}^m \frac {x_i}{1-e^{-x_i}}\right )\cdot ch(E).$$
Here $x_1,\ldots ,x_m$ are the roots of the almost complex structure.

For $\chi_y(M,E):=\sum _{p=0}^m \chi ^p(M,E)\cdot y^p\in \Z [y]$ one obtains by the index theorem
$$\chi _y(M,E)=Td_y(M,E):=\int _M \left (\prod _{i=1}^m \left ( (1+y\cdot e^{-x_i})\cdot \frac {x_i}{1-e^{-x_i}}\right )\right )\cdot ch (E). $$
Note that $Td_y(M,E)=\sum _p Td(M,\Lambda^p T^*M\otimes E)\cdot y^p$.

Now suppose a compact Lie group $G$ acts on $E$ and on $M$ preserving the almost complex structure. Then the Dolbeault complex refines to an equivariant complex and $\chi ^0(M,E)=Td(M,E)\in \Z $, resp. $\chi _y(M,E)=Td_y(M,E)\in \Z[y]$, refine to equivariant indices which we denote by $Td (M,E)_G\in R(G)$, resp. $Td_y (M,E)_G\in R(G)[y]$.

From now on we will assume that $c_1(M)\equiv 0 \bmod N$ for some fixed integer $N\geq 2$. Let $K=\det (T^*M)$ be the canonical line bundle. Since $c_1(M)\equiv 0 \bmod N$ there is a  complex line bundle $K^{1/N}$ on $M$ which represents an $N$th root of $K$. Now suppose the $G$-action lifts to $K^{1/N}$. Then one can consider the equivariant index $\chi_y(M,K^{\alpha /N})_G=Td_y(M,K^{\alpha /N})_G\in R(G)[y]$ for $\alpha\in \Z $.

For $G$ connected one has the following rigidity and vanishing results. If $\alpha =0$ then $Td_y(M,K^{\alpha /N})_G=Td_y(M)_G $ is rigid as shown by Lusztig and Kosniowski (see Theorem \ref{theoremToddrigid}).

If $0<\alpha <N $ and $y=0$ one has
$Td_y(M,K^{\alpha /N})_G=Td(M,K^{\alpha /N})_G=0$. This is due to Krichever \cite{Kr76} and Hattori \cite{Ha78}. Both statements, the rigidity statement, resp. the vanishing statement, can be shown by applying the cohomological form of Lefschetz fixed point formula \ref{fpfindK} to conclude that the equivariant index extends to a holomorphic function on the complex plane which is bounded, resp. which is bounded and vanishes at infinity.

Note that for $N=2$ and $\alpha =1$ one has $Td(M,K^{\alpha /N})=\hat A(M)$ and one recovers the Atiyah-Hirzebruch vanishing theorem for the $\hat A$-genus in the framework of almost complex manifolds.

The rigidity of $Td_y$ generalizes to the rigidity of elliptic genera of level $N$. This leads to higher vanishing results. We first give a definition of the elliptic genus of level $N$ following Witten's heuristic.

\subsection{${\boldsymbol {Td_y}}$ of the free loop space (elliptic genus of level ${\boldsymbol {N}}$)}

Let $M$ be an almost complex manifold with $c_1(M)\equiv 0 \bmod N$. As for the signature one can try to derive an invariant of $M$ starting with a hypothetical index $\chi _y$ on its free loop space ${\mathcal L}M$. Doing so one obtains (see \cite[\S 7.4, Appendix III]{HiBeJu92}) the series
$$Td_y({\mathcal L}M):=\sum _{n=0}^\infty Td_y(M,R_n) \cdot q^n\in \Z[y][[q]],$$
where $y=-e^{2\pi i/N}$ and the virtual bundles $R_n$ are determined by
$$\sum _{n=0}^\infty R_n \cdot q^n=\bigotimes_{n=1}^\infty \Lambda _{yq^n}T^*M \otimes \bigotimes_{n=1}^\infty \Lambda _{y^{-1}q^n}TM\otimes \bigotimes_{n=1}^\infty S_{q^n}
(TM+T^*M).$$
For example $R_0=1$ and $R_1=(1+y)\cdot T^*M +(1+y^{-1})\cdot TM$. Note that $Td_y({\mathcal L}M)$ can be computed from the Chern classes of $M$ and defines a genus for the unitary bordism ring. Also note that for $N=2$ the series $Td_y({\mathcal L}M)$ depends only on the Pontrjagin classes of $M$ and is equal to $sign({\mathcal {L}}M)$.

The series $Td_y({\mathcal L}M)$ converges for $\vert q\vert <1$ to a meromorphic function $\Phi_{1/N} (M)(\tau )$ on the upper half plane, where $q=e^{2\pi i \tau}$. Under the condition $y=-e^{2\pi i/N}$ it follows that $\Phi_{1/N} (M)$ is a modular function for the congruence subgroup
$$\Gamma _1(N):=\{A=(\begin{smallmatrix}a & b\\ c & d
\end{smallmatrix}) \in
SL_2(\Z)\mid (\begin{smallmatrix}a & b\\ c & d
\end{smallmatrix})\equiv (\begin{smallmatrix}1 & b\\ 0 & 1
\end{smallmatrix}) \bmod N\},$$
 (for an explanation of these facts see \cite{Wi86,Hi88,HiBeJu92}).

The genus $\Phi_{1/N}$ is called {\em elliptic genus of level $N$}. By construction $Td_y({\mathcal L}M)$ is the expansion of the modular function $\Phi_{1/N} (M)$ in the cusp $i\infty$.
Like for the elliptic genus of level $2$ the elliptic genus of level $N$ can be expanded in other cusps. Up to constants these expansions can be identified with
$$Td(M,\Lambda _y T^*M\otimes K^{\alpha /N}\otimes \bigotimes_{n=1}^\infty \Lambda _{yq^n}T^*M \otimes \bigotimes_{n=1}^\infty \Lambda _{y^{-1}q^n}TM\otimes \bigotimes_{n=1}^\infty S_{q^n}
(TM+T^*M))$$
for $ y=-e^{2\pi i \frac {\alpha \tau +\beta}N}$, $0\leq \alpha,\beta <N$, $(\alpha,\beta)\neq (0,0)$. The values of $\Phi_{1/N} (M)$ in the different cusps are given (up to constants) by
$$Td(M,K^{\alpha/N}),\; 0<\alpha <N \quad \text{and} \quad  Td_y(M),\; y=-e^{2\pi i\beta/N},\;  0< \beta <N .$$
Note that for $N=2$ these values are $Td(M,K^{1/2})=\hat A(M)$ and $Td_1(M)=sign(M)$.

\subsection{Rigidity of the elliptic genus of level $\boldsymbol N$}

Suppose a compact Lie group $G$ acts on $M$ preserving the almost complex structure and the $G$-action lifts to $K^{1/N}$. Then the elliptic genus of level $N$, $\Phi_{1/N} (M)$, refines to an equivariant elliptic genus $\Phi_{1/N} (M)_G$ and its expansions refine to series of virtual $G$-representations. For example the expansion $Td_y({\mathcal L}M)$ refines to $$ Td_y({\mathcal L}M)_G=\sum _{n=0}^\infty Td_y(M,R_n)_G \cdot q^n\in R(G)[y][[q]].$$

We are now in the position to state the
\begin{rigiditytheoremlevelN}[Witten, Hirzebruch, Bott-Taubes] If $G$ is connected then $\Phi_{1/N} (M)_G$ and its expansions are rigid, i.e. each coefficient in the expansions is constant as a character of $G$.\proofend\end{rigiditytheoremlevelN}

The rigidity of the elliptic genus of level $N$ was conjectured by Witten \cite{Wi86} and then shown by Hirzebruch \cite{Hi88} and Bott-Taubes \cite{BoTa89} (see also \cite[Appendix III]{HiBeJu92}). The proof uses the fixed point formula and elliptic function theory.

\subsection{Higher order vanishing for level ${\boldsymbol {N}}$}
Let $M$ be as before an almost complex manifold with $c_1(M)\equiv 0 \bmod N$. Suppose $S^1$ acts on $M$ preserving the almost complex structure. If the $S^1$-action does not lift to $K^{1/N}$ then a connected covering of $S^1$ will lift. In this case one can show, using the Lefschetz fixed point formula and the rigidity theorem, that the equivariant elliptic genus of level $N$ vanishes identically (see \cite[Thm. on page 58]{Hi88}). Hence, for the following vanishing results we may, without mentioning, always assume that the $S^1$-action lifts to $K^{1/N}$.

As pointed out by Hirzebruch in \cite[\S 11]{Hi88} the rigidity of the elliptic genus of level $N$ can be used to show the vanishing of $Td(M,K^{\alpha/N}),\; 0<\alpha <N$. This is another manifestation of the {\em rigidity implies vanishing} phenomenon and provides an alternative proof of the theorem of Krichever \cite{Kr76} and Hattori \cite{Ha78}. We note that this kind of reasoning is completely analogous to the one for the vanishing of the $\hat A$-genus using the rigidity of the elliptic genus of level $2$.

Recall from Theorems \ref{HiSlvanishing} and \ref{highervanishing} that for certain $S^1$-actions next to the $\hat A$-genus additional coefficients in the expansion of the elliptic genus of level $2$ vanish. It turns out that the idea of the proof of these higher vanishing results which is given in \cite[\S 5]{De07} also applies to elliptic genera of level $N$. This leads to new vanishing results for $\Phi_{1/N} (M)$ and applications to positively curved manifolds (see the forthcoming preprint \cite{De} for details). We close this section with two sample applications.

\begin{theorem} Let $M$ be an almost complex manifold with $c_1(M) \equiv 0 \bmod N$. Suppose $S^1$ acts on $M$ preserving the almost complex structure and let $\sigma \in S^1$ be the element of order $2$. If $\codim M^{\sigma}\geq 6$ then $Td(M,T^*M\otimes K^{\alpha/N})$ vanishes for $0<\alpha < N/2$.\proofend\end{theorem}

\begin{theorem} Let $M$ be a connected almost complex manifold with $c_1(M) \equiv 0 \bmod N$. Suppose $M$ has positive sectional curvature and a torus $T$ of rank $3$ acts isometrically and effectively on $M$ preserving the almost complex structure. Then  $Td(M,T^*M\otimes K^{\alpha/N})$ vanishes for $0<\alpha < N/2$.\proofend\end{theorem}
Similar results hold if one replaces positive sectional by positive $k$th Ricci curvature.

\begin{finalremarks}\begin{enumerate}
\item If $c_1(M)=0$ the rigidity and vanishing results above can be used to prove corresponding statements for the $SU$-elliptic genus.
\item Based on \cite{De00} similar results can be shown for spin$^c$ manifolds.
\item If $M$ is a Fano manifold with $c_1(M) \equiv 0 \bmod N$ then some of the higher vanishing results for the elliptic genus of level $N$ hold without any symmetry assumption.
\end{enumerate}
\end{finalremarks}

\end{document}